\documentclass[11pt]{article}
\usepackage{latexsym,amssymb,amsmath,a4wide,theorem}
\usepackage{graphicx}
\usepackage{subfig}
\usepackage{caption}
\usepackage{color}
\usepackage{extarrows}
\usepackage{esint}

\numberwithin{equation}{section}
\numberwithin{figure}{section}
\newtheorem{theorem}{Theorem}[section]
\newenvironment{proof}{{ \it Proof:\quad}}{\hfill $\blacksquare$\par}
\newtheorem{lemma}{Lemma}[section]
\newtheorem{remark}{Remark}[section]
\newtheorem{proposition}{Proposition}[section]

\begin{document}
\title{\bf Sharp Interface Limit for  Compressible Immiscible Two-Phase Dynamics with Relaxation}
\author{\large Yazhou Chen$^{1}$, Yi Peng$^{1}$, Qiaolin He$^{2}$,  Xiaoding Shi$^{1*}$\\
\scriptsize$^{1}$ {College of Mathematics and Physics, Beijing University of
Chemical Technology, Beijing 100029, China}\\
\scriptsize$^{2}$ {School of Mathematics, Sichuan University, Chengdu, Sichuan, 610065,  China}
}

\date{}
\maketitle

\begin{abstract}
In this paper, the  compressible immiscible two-phase flow with relaxation is investigated, this model can be regarded as a natural modification of Jin-Xin relaxation scheme proposed and developed by S.Jin and Z.P.Xin([Comm.Pure Appl.Math., 48,1995]) in
view of the numerical approximation of conservation laws. Given any entropy solution consists of two different families of shocks interacting at some positive time for the standard two-phase compressible Euler equations, it is proved that such entropy solution is the  sharp interface limit for a family global strong solutions of the modified Jin-Xin relaxation scheme for Navier-Stokes/Allen-Cahn system,  here the relaxation time is selected as the thickness of the interface,  weighted estimation and improved antiderivative method are used in the proof. Moreover, the simulation results are given by this modified Jin-Xin relaxation scheme method. Both numerical and theoretical results show that, the interacting shock waves can pass through the interface without any effect.

\end{abstract}

\

\noindent{\bf Keywords}: Sharp Interface Limit, Compressible Immiscible Two-Phase Dynamics, Shock Wave, Rarefaction Wave, Jin-Xin Relaxation.

\

\noindent{\textbf{MSC}}: 35Q35, 35B65, 76N10, 35M10, 35B40, 35C20, 76T30

\renewcommand{\thefootnote}{\fnsymbol{footnote}}
\footnotetext[1]{Corresponding author: shixd@mail.buct.edu.cn, (Xiaoding Shi)}

\section{\normalsize Introduction and Main Results}

\ \ \ \
The interfacial dynamic for the immiscible two-phase flow has  attracted extensive research for hundreds years. The study of interface structure, phase behavior, and dynamics is fascinating. In this paper, the following Navier-Stokes/Allen-Cahn system modeling the compressible immiscible two-phase flow is considered(see Blesgen \cite{B1999},  Heida-M$\mathrm{\acute{a}}$lek-Rajagopal  \cite{HMR2012} and the references therein):
\begin{equation}\label{original NSAC}
\left\{\begin{array}{llll}
\displaystyle \rho_{t}+\textrm{div}(\rho \mathbf{u})=0,\\
\displaystyle (\rho \mathbf{u})_{t}+\mathrm{div}\big(\rho \mathbf{u}\otimes \mathbf{u}\big)=-\mathrm{div}\Big(\mathbb{P}+\rho\nabla\chi\otimes\frac{\partial f}{\partial\nabla\chi}\Big),
  \\
\displaystyle(\rho\chi)_{t}+\mathrm{div}\big(\rho\chi \mathbf{u}\big)=-L_d(\epsilon)\mu,\quad
\displaystyle\rho\mu=\rho\frac{\partial f}{\partial \chi}-\mathrm{div}\big(\rho\frac{\partial f}{\partial \nabla\chi}\big),
\end{array}\right.
\end{equation}
where $\rho$, $\mathbf{u}$, $\chi$  are represented separately the total density, the  average velocity and the phase field. $\mu(\mathbf{x},t)$  is the chemical potential,  $f$  the  phase-phase interfacial free energy density, $\mathbb{P}$ the  stress tensor,  they have the following expression (see Heida-M$\mathrm{\acute{a}}$lek-Rajagopal,  \cite{HMR2012}Lowengrub-Truskinovsky \cite{LT-1998}):
\begin{align}\label{free energy density and total energy density}
 &f(\rho,\chi,\nabla\chi)\overset{\text{def}}{=}\frac{1}{4}(1-\chi^2)^2+\frac{\epsilon^2}{2\rho}|\nabla \chi|^2,\ \ \ \mathbb{P}=\Big(R\rho^{\gamma}-\frac{\epsilon^2}2|\nabla\chi|^2\Big)\mathbb{I},
 \end{align}
where $R$ is the gas constant, and $\gamma>1$ the adiabatic exponent. Without loss of generality, we assume $R=1$. $\epsilon>0$  represents the interfacial thickness of immiscible two-phase flow, and $L_d(\epsilon)=O(\frac{1}{\epsilon})$  the phenomenological mobility coefficient related to the speed at which the system  approaches an equilibrium configuration, which means that the diffusivity $L_d(\epsilon)$ decreases as the interfacial thickness $\epsilon$  increases.

The Navier-Stokes/Allen-Cahn system \eqref{original NSAC}  is widely used to describe the compressible immiscible two-phase flow. In particular, this system can well describe the phase transition and phase separation of immiscible two-phase flow. Nevertheless, since the physical scale of the thickness for interface of the two-phase flow is actually quite small, the scale is so small that even computation cannot select such small physical scales, so the sharp interface limit of this system \eqref{original NSAC} becomes an important and challenging problem. Compared with many researches on the well-posedness of this system, (see \cite{CHHS2021}-\cite{FPRS2010}, \cite{K2012},\cite{LYZ-2018}-\cite{LYZ-2020},\cite{YZ-2018}, etc.), not so much results on the problem of  sharp interface limit, even for the computational simulation,  there are only a few  related works, for example,  Witterstein \cite{W-2010} used the asymptotic matching method to analyze this problem, and showed that the sharp-interface limit is the free boundary problem for two-phase compressible Navier-Stokes equations formally.

The motivation of this paper is to construct a computational scheme for inviscid compressible Navier-Stokes/Allen-Cahn \eqref{original NSAC}, and the convergence of the scheme is also analyzed. Thanks for the Jin-Xin relaxation schemes \cite{JinXin1995}, which proposed for the numerical approximation of  conservation laws, a modified  Jin-Xin relaxation schemes \eqref{Jin-Xin relaxation system} is constructed for the system  \eqref{original NSAC} in this paper. In this modified relaxation scheme we constructed, The thickness of the interface $\epsilon$ is selected as the relaxation parameter, so that the sharp interface limit is transformed into the relaxation limit problem.  As well as we known, for the smooth solution, the relaxation limit problem can be studied by scaling method. However,  the solution for \eqref{Jin-Xin relaxation system} in general develops singularity in the finite time such as shock waves. Thus it is important but difficult to study the singular limit as $\epsilon\rightarrow0^+$  in the case when the solution of system \eqref{Jin-Xin relaxation system} does not possess smooth solution.
There are a lot of results on the analysis of the large time behavior for the solution to the Jin-Xin scheme of the conservation laws, one may refer to \cite{HWWY}, \cite{SZ}--\cite{SYZ}, \cite{YZ00}-\cite{Z02} and the references therein.

As a beginning to study the sharp interface limit \eqref{Jin-Xin relaxation system}, we first to investigate the  asymptotic nonlinear stability of a superposition of interacting relaxation shock waves belong to different of families for \eqref{Jin-Xin relaxation system}.  We want to know that if the interacting relaxation shock waves crosses the interface of immiscible two-phase flow, whether the interface affects the properties of the shock waves. In particular, the sharp interface limit to the modified Jin-Xin scheme we constructed here.

\vspace{0.2cm}\noindent\textbf{Notation.} Throughout this paper,
$L^2(\mathbb{R})$, $W_2^k(\mathbb{R})$ represent the usual Sobolev spaces on $\mathbb{R}$, with norms  $\|\cdot\|$ and $\|\cdot\|_{k}$, respectively.   $C, c$ used to denote the constants which are
independent of $x,t$.   We will employ the notation $a\lesssim b$ to mean that $a\leq  Cb$ for a universal constant $C>0$ that only depends on the parameters coming from the problem.

\vspace{0.2cm}

Letting  $v=\frac1\rho$ the specific volume, under the Lagrange coordinate system,  system \eqref{original NSAC} in one-dimensional can be rewritten as:
\begin{equation}\label{NSFAC-Lagrange}
\left\{\begin{array}{llll}
\displaystyle v_t-u_x=0,\\
\displaystyle u_t+p_x=-\frac{\epsilon^2}{2}(\frac{\chi_x^2}{v^2})_x, \\
\displaystyle \chi_t=-\frac{v}{\epsilon}\mu,\ \ \ \displaystyle \mu=(\chi^3-\chi)-\epsilon^2 \Big(\frac{\chi_x}{v}\Big)_x,
\end{array}\right.
\end{equation}
where $  p=\frac{1}{v^{\gamma}}$, setting
\begin{equation}\label{definition of U and F}
  \mathbf{U}=\left(
               \begin{array}{c}
                 v \\
                 u
               \end{array}
             \right),\quad
              \mathbf{F}=\left(
    \begin{array}{c}
      -u \\
   \displaystyle   p     \end{array}
  \right),
\end{equation}
formally, as $\epsilon\rightarrow0^+$, the system \eqref{NSFAC-Lagrange} tends to the following standard two-phase compressible inviscid Navier-Stokes equations:
\begin{equation}\label{Euler Equations for two-phase flow}
\left\{\begin{array}{llll}
\displaystyle \mathbf{U}_t+\mathbf{F}(\mathbf{U})_x=0,\ & \mathrm{in}\ \Omega^{\pm}(t),\\
\displaystyle \chi=\pm1, & \mathrm{in}\ \Omega^{\pm}(t),
\end{array}\right.
\end{equation}
where $ \Omega^{-}({t})=\big\{{x}\in\mathbb{R}\big|\chi({x},t)=-1\big\}$, $\Omega^+({t})=\big\{x\in\mathbb{R}\big|\chi({x},{t})=1\big\}$,
 and $\Gamma({t})=  \mathbb{R}\backslash\big\{\Omega^-(t)\cup\Omega^+(t)\big\}$, $\mathrm{meas}\Gamma(t)=0$, $t\geq0$.

In the following, we consider the interaction of the two shock waves belong to different family for  the system \eqref{Euler Equations for two-phase flow} with respect to the following initial condition:

\begin{eqnarray} \label{initial for Euler}
\mathbf{U}(x,0)=\left(
               \begin{array}{c}
                 v_0 \\
                 u_0
               \end{array}
             \right)=
\left\{ \begin{aligned}
& \mathbf{U}_-=(v_-,u_-)^\mathbb{T},\,\,x<0,  \\
& \mathbf{U}^*=(v^*,u^*)^\mathbb{T}, \,\,0<x<1, \\
& \mathbf{U}_+=(v_+,u_+)^\mathbb{T}, \,\,x>1.
\end{aligned} \right.
\end{eqnarray}
where $\mathbb{T}$ denotes the transpose and all vectors are column ones. Noting that the eigenvalues of the  system \eqref{Euler Equations for two-phase flow} are
\begin{eqnarray}
\lambda_1=-\sqrt{-p'_v}<0,\ \ \lambda_2=\sqrt{-p'_v}>0,
\end{eqnarray}
we assume that, $v_{\pm}>0$, $(v_\pm,u_\pm)$, $(v^*,u^*)$ and $s_1,s_2$ are constants satisfying Rankine-Hugoniot conditions
\begin{eqnarray}\label{R-H Before the collision}
&-s_1(\mathbf{U}_+-\mathbf{U}^*)+\mathbf{F}(\mathbf{U}_+)-\mathbf{F}(\mathbf{U}^*)=0, \\
& -s_2(\mathbf{U}^*-\mathbf{U}_-)+\mathbf{F}(\mathbf{U}^*)-\mathbf{F}(\mathbf{U}_-)=0,
\end{eqnarray}
and the Lax's shock condition given below:
\begin{eqnarray}\label{entopy condition 1}
\lambda_1(\mathbf{U}_+)<s_1<\lambda_1(\mathbf{U}^*)<0<\lambda_2(\mathbf{U}^*)<s_2<\lambda_2(\mathbf{U}_-).
\end{eqnarray}
By simple calculation, we get the entropy solution $(\mathcal{V},\mathcal{U})$ expression of system \eqref{Euler Equations for two-phase flow},\eqref{initial for Euler} as following,
 \begin{eqnarray}\label{solution before t0}
(\mathcal{V},\mathcal{U})(x,t)\xlongequal{t\leq t_0}
\left\{ \begin{aligned}
& \mathbf{U}_-,\quad\,\,x<s_2t,  \\
& \mathbf{U}^*, \quad\,\,s_2t<x<s_1t+1,  \\
& \mathbf{U}_+, \quad\,\,x>s_1t+1,
\end{aligned} \right.
\end{eqnarray}
and
\begin{eqnarray} \label{solution after t0}
(\mathcal{V},\mathcal{U})(x,t)\xlongequal{t>t_0}
\left\{ \begin{aligned}
& \mathbf{U}_-=(v_-,u_-)^{\mathbb{T}},\,\quad\,x-x_0<\tilde{s}_1(t-t_0),   \\
& \mathbf{\tilde{U}}^*=(\tilde{v}^*,\tilde{u}^*)^{\mathbb{T}},\quad\,\,\,\,\tilde{s}_1(t-t_0)<x-x_0<\tilde{s}_2(t-t_0),\\
& \mathbf{U}_+=(v_+,u_+)^{\mathbb{T}}, \quad\,\,\,x-x_0>\tilde{s}_2(t-t_0),
\end{aligned} \right.
\end{eqnarray}
where $Q\triangleq(x_0,t_0)=(\frac{s_2}{s_2-s_1},\frac1{s_2-s_1})$ is the intersection of the two incoming reverse shock waves $S_1,S_2$.  Obviously, the physical meaning of the solution \eqref{solution before t0} is,
two incoming reverse shock waves $S_1,S_2$ are initially given,  $S_1$ represents the 1-shock which connects $\mathbf{U}^*$ as the left state
and $\mathbf{U}_+$ as the right state with traveling wave velocity $s_1<0$,  $S_2$  the 2-shock which connects $\mathbf{U}_-$ as the left state
 and $\mathbf{U}^*$ as the right state with traveling wave velocity $s_2>0$, and these two shocks must interact at point $Q$.
After the interaction of two reverse shock waves, two outgoing shocks $\tilde{S}_1$,$\tilde{S}_2$  are generated. Assume that the intermediate state between $\tilde{S}_1$ and $\tilde{S}_2$ is $\tilde{\mathbf{U}}^*$,  these two intermediate states are determined by the Rankine-Hugoniot condition as following
\begin{eqnarray}
&  -\tilde{s}_1(\mathbf{\tilde{U}}^*-\mathbf{U}_-)+\mathbf{F}(\mathbf{\tilde{U}}^*)-\mathbf{F}(\mathbf{U}_-)=0,  \\
& -\tilde{s}_2(\mathbf{U}_+-\mathbf{\tilde{U}}^*)+\mathbf{F}(\mathbf{U}_+)-\mathbf{F}(\mathbf{\tilde{U}}^*)=0,
\end{eqnarray}
and the  Lax's shock condition given below
\begin{eqnarray}
\lambda_1(\mathbf{\tilde{U}}^*)<\tilde{s}_1<\lambda_1(\mathbf{U}_-)<0<\lambda_2(\mathbf{U}_+)<\tilde{s}_2<\lambda_2(\mathbf{\tilde{U}}^*).
\end{eqnarray}

It is well known that the system \eqref{Euler Equations for two-phase flow} is a nonlinear hyperbolic partial differential equations  with free boundaries, \eqref{NSFAC-Lagrange}  can be regarded as an approximation of the problem  \eqref{Euler Equations for two-phase flow} in a sense.
 So the singular limit problem of the interface thickness approaching zero becomes particularly important both computationally and analytically. In order to solve this sharp interface limit, in this paper, we give an relaxation scheme for calculation and analysis based on the famous Jin-Xin relaxation method \cite{JinXin1995}. We choose the thickness of the interface as the relaxation  factor, then the modified Jin-Xin relaxation scheme for \eqref{NSFAC-Lagrange} is constructed as follows:
\begin{equation}\label{Jin-Xin relaxation system}
\left\{\begin{array}{llll}
\displaystyle \mathbf{U}_t+\mathbf{V}_x=\left(
               \begin{array}{c}
                 0 \\
                 -\frac{\epsilon^2}{2}(\frac{\chi_x^2}{v^2})_x
                                \end{array}
             \right),
             \\
\displaystyle \mathbf{V}_t+a^2\mathbf{U}_{x}=\frac{1}{\epsilon}[\mathbf{F}(\mathbf{U})-\mathbf{V}],\\
\displaystyle \chi_t=-\frac{v}{\epsilon}\mu,\ \ \ \displaystyle \mu=(\chi^3-\chi)-\epsilon^2\Big(\frac{\chi_x}{v}\Big)_x,
\end{array}\right.
\end{equation}
where, the positive constant $a$  should be chosen  that satisfies the following sub-characteristic condition (see \cite{JinXin1995}, \cite{LX97}):
\begin{equation}\label{sub-characteristic condition}
 \sqrt{-p'(v)}<a.
\end{equation}
In order to analyze the singularity limit $\epsilon\rightarrow0$ for the problem \eqref{Jin-Xin relaxation system}, the scaling method is adopted as following:
\begin{equation}\label{Scaling transformation}
y=\frac{x-x_0}{\epsilon},\qquad \tau=\frac{t-t_0}{\epsilon},
\end{equation}
then \eqref{Jin-Xin relaxation system} can be rewritten as
\begin{equation}\label{Jin-Xin relaxation system y}
\left\{\begin{array}{llll}
\displaystyle \mathbf{U}_\tau+\mathbf{V}_y=\mathbf{H}_y,\ \ \\
\displaystyle \mathbf{V}_\tau+a^2\mathbf{U}_{y}=[\mathbf{F}(\mathbf{U})-\mathbf{V}],\\
\displaystyle \chi_\tau=-v\mu,\ \ \ \displaystyle \mu=(\chi^3-\chi)-\Big(\frac{\chi_x}{v}\Big)_x,
\end{array}\right.
\end{equation}
where
\begin{equation}\label{H}
 \mathbf{H}(y,\tau)=\Big(0,-\frac{1}{2}\big(\frac{\chi_y^2}{v^2}\big)\Big)^\mathbb{T}.
\end{equation}
Differentiating  \eqref{Jin-Xin relaxation system y}$_2$ with respect to  $y$,  using \eqref{Jin-Xin relaxation system y}$_1$, one has
\begin{equation}\label{Relaxation system}
\left\{\begin{array}{llll}
\displaystyle \mathbf{U}_\tau+\mathbf{F}_y(\mathbf{U})=(a^2 \mathbf{U}_{yy}-\mathbf{U}_{\tau\tau})+\mathbf{H}_{y\tau}+\mathbf{H}_y, \\
\displaystyle \chi_\tau=-v\mu,\ \ \ \displaystyle \mu=(\chi^3-\chi)-\Big(\frac{\chi_x}{v}\Big)_x,
\end{array}\right.
\end{equation}
From this, we associate to construct the relaxation shock solutions $\mathbf{U}(\xi)=\mathbf{U}(y-s\tau)$ of the following wave equations for the shock waves $(\mathcal{V},\mathcal{U})$ of system \eqref{Euler Equations for two-phase flow},\eqref{initial for Euler},
\begin{equation}\label{relaxation shock wave}
 \mathbf{U}_\tau+\mathbf{F}_y(\mathbf{U})=(a^2 \mathbf{U}_{yy}-\mathbf{U}_{\tau\tau}).
\end{equation}
Letting $\mathbf{r}_i(\mathbf{U})$, $\mathbf{l}_i(\mathbf{U})$($i=1, 2$) be the right and left eigenvectors of $\mathbf{F}'(\mathbf{U})$ such that
\begin{eqnarray}
\mathbf{F}'(\mathbf{U})\mathbf{r}_i(\mathbf{U})=\lambda_i(\mathbf{U})\mathbf{r}_i(\mathbf{U}), \quad \mathbf{l}_i(\mathbf{U})\mathbf{F}'(\mathbf{U})=\lambda_i(\mathbf{U})\mathbf{l}_i(\mathbf{U}),\,\,\,\mathrm{with}\,\,\,i=1,2,
\end{eqnarray}
where $\mathbf{r}_j(\mathbf{U})$ and $\mathbf{l}_i(\mathbf{U})$ satisfy that
\begin{eqnarray}
\mathbf{l}_i(\mathbf{U})\mathbf{r}_j(\mathbf{U})=\delta_{ij}, \,\,\,\mathrm{with}\,\,\,i, j=1,2.
\end{eqnarray}
Constructing the following matrix $\mathbf{L}$, $\mathbf{R}$, $\mathbf{\Lambda}$:
\begin{eqnarray}
\mathbf{L}=(\mathbf{l}_1,\mathbf{l}_2)^\mathbb{T},\,\,\,\mathbf{R}=(\mathbf{r}_1,\mathbf{r}_2),\,\,\,\mathbf{\Lambda}=\mathrm{diag}(\lambda_1,\lambda_2),
\end{eqnarray}
and then one has
\begin{eqnarray}
\mathbf{L}\mathbf{F}'(\mathbf{U})\mathbf{R}=\mathbf{\Lambda},\,\,\,\mathbf{L}\mathbf{R}=\mathbf{I}.
\end{eqnarray}
After a simple calculation, one gets the 1-relaxation shock wave $\mathbf{U}^{S_{1}}=(V^{S_1},U^{S_1})(y-s_1\tau)$ exists uniquely up to a shift, which connecting $\mathbf{U}^*$ on the left and $\mathbf{U}^+$ on the right, and satisfies
\begin{eqnarray} \label{1-relaxation shock wave}
\left\{ \begin{aligned}
&-s_1\mathbf{U}^{S_1}_{\xi_1}+\mathbf{F}(\mathbf{U}^{S_1})_{\xi_1}=(a^2-s_1^2) \mathbf{U}^{S_1}_{\xi_1\xi_1},\\
&\mathbf{U}^{S_1}(-\infty)=\mathbf{U}^*,\ \ \mathbf{U}^{S_1}(+\infty)=\mathbf{U}_+,
\end{aligned} \right.
\end{eqnarray}
where  $\xi_1=y-s_1\tau$. To fix the relaxation shock wave, after calculation, we construct the following initial condition
\begin{equation} \label{initial of v S1}
\mathbf{U}^{S_1}(0)=\frac{1}{2}(\mathbf{U}_++\mathbf{U}^*),
\end{equation}
and for the 2-relaxation shock wave $\mathbf{U}^{S_{2}}=(V^{S_2},U^{S_2})(y-s_2\tau)$, which connecting $\mathbf{U}_-$ on the left and $\mathbf{U}^*$ on the right, one obtains
\begin{eqnarray} \label{2-relaxation shock wave}
\left\{ \begin{aligned}
& -s_2\mathbf{U}^{S_2}_{\xi_2}+\mathbf{F}(\mathbf{U}^{S_2})_{\xi_2}=(a^2-s_2^2) \mathbf{U}^{S_2}_{\xi_2\xi_2},\\
&\mathbf{U}^{S_2}(-\infty)=\mathbf{U}_-,\quad\mathbf{U}^{S_2}(+\infty)=\mathbf{U}^*,\\
&\mathbf{U}^{S_2}(0)=\frac{1}{2}(\mathbf{U}^*+\mathbf{U}_-),
\end{aligned} \right.
\end{eqnarray}
where  $\xi_2=y-s_2\tau$.
 Similarly,  the 1-relaxation outgoing shock $\mathbf{\tilde{U}}^{\tilde{S}_1}=(\tilde{V}^{\tilde{S}_1},\tilde{U}^{\tilde{S}_1})(y-\tilde{s}_1\tau)$, which connecting $\mathbf{U}_-$ on the left and $\tilde{\mathbf{U}}^*$ on the right, and the 2-relaxation outgoing shock $\mathbf{\tilde{U}}^{\tilde{S}_2}=(\tilde{V}^{\tilde{S}_2},\tilde{U}^{\tilde{S}_2})(y-\tilde{s}_2\tau)$, which connecting $\mathbf{U}^*$ and $\mathbf{U}_+$, are constructed as following
 \begin{eqnarray}
\left\{ \begin{aligned}
& -\tilde{s}_1\mathbf{\tilde{U}}^{\tilde{S}_1}_{\tilde{\xi}_1}+\mathbf{F}(\mathbf{\tilde{U}}^{\tilde{S}_1})_{\tilde{\xi}_1}=(a^2-\tilde{s}_1^2) \mathbf{\tilde{U}}^{\tilde{S}_1}_{\tilde{\xi}_1\tilde{\xi}_1},\\
&\mathbf{\tilde{U}}^{\tilde{S}_1}(-\infty)=\mathbf{U}_-,\quad\mathbf{\tilde{U}}^{\tilde{S}_1}(+\infty)=\tilde{\mathbf{U}}^*,\\
&\mathbf{\tilde{U}}^{\tilde{S}_1}(0)=\frac{1}{2}(\mathbf{U}_-+\mathbf{U}_-^*),
\end{aligned} \right.
\end{eqnarray}
and
\begin{eqnarray}
\left\{ \begin{aligned}
& -\tilde{s}_2\mathbf{\tilde{U}}^{\tilde{S}_2}_{\tilde{\xi}_2}+\mathbf{F}(\mathbf{\tilde{U}}^{\tilde{S}_2})_{\tilde{\xi}_2}=(a^2-\tilde{s}_2^2) \mathbf{\tilde{U}}^{\tilde{S}_2}_{\tilde{\xi}_2\tilde{\xi}_2},\\
&\mathbf{\tilde{U}}^{\tilde{S}_2}(-\infty)=\mathbf{U}^*,\quad\mathbf{\tilde{U}}^{\tilde{S}_2}(+\infty)=\mathbf{U}_+,\\
&\mathbf{\tilde{U}}^{\tilde{S}_2}(0)=\frac{1}{2}(\mathbf{U}_+^*+\mathbf{U}_+).
\end{aligned} \right.
\end{eqnarray}
Moreover, the strengths of the two incoming shocks are defined as follows
\begin{equation}
\delta_1=|(v_+-v^*,u_+-u^*)| ,\quad\delta_2=|(v^*-v_-,u^*-u_-)|,
\end{equation}
and setting
\begin{equation}\label{The strength of a wave}
  \delta=\min\{\delta_1,\delta_2\}.
\end{equation}
 In the same way, the strengths of the outgoing waves are defined as
\begin{equation}\label{The intensity of the reflected wave}
\begin{split}
&\tilde{\delta}_1=|(\tilde{v}^*-v_-,\tilde{u}^*-u_-)|, \ \
\tilde{\delta}_2=|(v_+-\tilde{v}^*,u_+-\tilde{u}^*)|,
\end{split}
\end{equation}
and satisfying(see \cite{S94})
\begin{eqnarray}\label{A change in the intensity of a wave}
\tilde{\delta}_1=\delta_1+O(1)\delta_1\delta_2, \quad\tilde{\delta}_2=\delta_2+O(1)\delta_1\delta_2,
\end{eqnarray}
for $\delta$ small enough. In the following, the decay property of relaxed shock wave and the properties of eigenvalues are given as following, which can be found in \cite{HWWY}, \cite{Ze09}  and the references therein.

\begin{lemma} \label{estimate of viscous shock}
There are positive constant $C_0$ and $c_0$, for $i=1,2$ and $\forall \tau\geq-\frac{t_0}{\epsilon}$ such that
\begin{eqnarray}
&&|\mathbf{U}^{S_1}-\mathbf{U}^*|\lesssim\delta_1e^{-c_0\delta_1|y-s_1\tau|},y\leq s_1\tau; \qquad|\mathbf{U}^{S_1}-\mathbf{U}_+|\lesssim\delta_1e^{-c_0\delta_1|y-s_1\tau|},y\geq s_1\tau; \notag \\
&&|\mathbf{U}^{S_2}-\mathbf{U}_-|\lesssim\delta_2e^{-c_0\delta_2|y-s_2\tau|},y\leq s_2\tau; \qquad|\mathbf{U}^{S_2}-\mathbf{U}^*|\lesssim\delta_2e^{-c_0\delta_2|y-s_2\tau|},y\geq s_2\tau ; \notag \\
&&|\mathbf{\tilde{U}}^{\tilde{S}_1}-\mathbf{U}_-|\lesssim\tilde{\delta}_1e^{-c_0\tilde{\delta}_1|y-\tilde{s}_1\tau|},y\leq \tilde{s}_1\tau;\qquad |\mathbf{\tilde{U}}^{\tilde{S}_1}-\mathbf{\tilde{U}}^*|\lesssim\tilde{\delta}_1e^{-c_0\tilde{\delta}_1|y-\tilde{s}_1\tau|},y\geq \tilde{s}_1\tau;    \notag \\
&&|\mathbf{\tilde{U}}^{\tilde{S}_2}-\mathbf{\tilde{U}}^*|\lesssim\tilde{\delta}_2e^{-c_0\tilde{\delta}_2|y-\tilde{s}_2\tau|},x\leq \tilde{s}_2\tau;\qquad |\mathbf{\tilde{U}}^{\tilde{S}_2}-\mathbf{U}_+|\lesssim\tilde{\delta}_2e^{-c_0\tilde{\delta}_2|y-\tilde{s}_3\tau|},y\geq \tilde{s}_2\tau; \notag \\
&&|\mathbf{U}^{S_i}_y|\lesssim\delta_i^2e^{-c_0\delta_i|\xi_i|},\quad|\mathbf{U}^{\tilde{S}_i}_y|\leq C\tilde{\delta}_i^2e^{-c_0\delta_1|\tilde{\xi}_i\tau|}; \quad\frac d{d\xi_i}\lambda_i(\mathbf{U}^{S_i})<0, \quad \frac d{d\tilde{\xi}_i}\lambda_i(\mathbf{\tilde{U}}^{\tilde{S}_i})<0;\notag\\
&&\qquad\notag\\
&&\big|\frac{d}{d\xi_i}\mathbf{U}^{S_i}\big|\lesssim\big|\frac{d}{d\xi_i}\lambda_i(\mathbf{U}^{S_i})\big|\lesssim\big|\frac{d}{d\xi_i}\mathbf{U}^{S_i}\big|;\ \big|\frac{d}{d\tilde{\xi}_i}\mathbf{\tilde{U}}^{\tilde{S}_i}\big|\lesssim\big|\frac{d}{d\tilde{\xi}_i}\lambda_i(\mathbf{\tilde{U}}^{\tilde{S}_i})\big|\lesssim\big|\frac{d}{d\tilde{\xi}_i}\mathbf{\tilde{U}}^{\tilde{S}_i}\big|;\notag\\
&&\int_{\mathbb{R}}\big|\frac{d}{d\xi_i}\lambda_i(\mathbf{U}^{S_i})\big|d\xi_i=O(\delta_i),\qquad \int_{\mathbb{R}}\big|\frac{d}{d\tilde{\xi}_i}\lambda_i(\mathbf{\tilde{U}}^{\tilde{S}_i})\big|d\tilde{\xi}_i=O(\tilde{\delta}_i);\notag\\
&&\big|\frac{d}{d\xi_i}\lambda_i(\mathbf{U}^{S_i})\big|\lesssim\delta_i^2e^{-C\delta_i|\xi_i|},\qquad \big|\frac{d}{d\tilde{\xi}_i}\lambda_i(\mathbf{\tilde{U}}^{\tilde{S}_i})\big|\lesssim\tilde{\delta}_i^2e^{-C\tilde{\delta}_i|\tilde{\xi}_i|}.\notag
\end{eqnarray}
\end{lemma}

\begin{theorem}(Main Result)  \label{main thm}
Assume that $(\mathcal{V},\mathcal{U})$ is the entropy solution of system \eqref{Euler Equations for two-phase flow}-\eqref{initial for Euler}  given by (\ref{solution before t0})-(\ref{solution after t0}). Under the sub-characteristic condition (\ref{sub-characteristic condition}),  there exist  positive constants $\delta_0$ and $M_0$ such that,
if the wave strength satisfies $\delta\leq\delta_0$ and the initial  perturbation $\|\chi_0^2-1\|\leq M_0$, the system (\ref{Jin-Xin relaxation system}) with the initial data \eqref{relaxation system before t0}$_3$ admits a family of global
smooth solution $(v,u,\chi)$. Moreover,  before the interaction time $t=t_0$, it holds that
\begin{equation}\label{convegence rate}
\lim_{\epsilon\rightarrow0}\|v^{\epsilon}-\mathcal{V},u^{\epsilon}-\mathcal{U},\chi^{2}-1\|_{L^\infty(\Sigma_h)}=0,
\end{equation}
where $\Sigma_h=\{(x,t)\mid  |x-s_2t|\geq h, \, |x-s_1t-1|\geq h,\, 0\leq t\leq t_0-h\}$ for any positive constant $h>0$.
After the interaction time $t_0$, it holds that
\begin{equation}\label{stability}
\lim_{\epsilon\rightarrow0}\|v-\mathcal{V},u-\mathcal{U},\chi^{2}-1\|_{L^\infty(\widetilde{\Sigma}_{\tilde{h}})}=0,
\end{equation}
where
$\widetilde{\Sigma}_{\tilde{h}}=\{(x,t)\mid |(x-x_0)-\tilde{s}_2(t-t_0)|\geq \tilde h, \,  |x-\tilde{s}_1(t-t_0)|\geq \tilde{h},\,  t_0
+\tilde {h}\leq t<+\infty\}$
for any positive constant $\widetilde h>0$.
\end{theorem}

\begin{remark}
The system \eqref{Jin-Xin relaxation system} can be seen as a modified Jin-Xin relaxation scheme which proposed by Jin-Xin in \cite{JinXin1995} for hyperbolic conservation law.
\end{remark}

\begin{remark}
Given any entropy solution of compressible immiscible inviscid two-phase flow which consists of two different families of shocks interacting at some positive time,
 Theorem 1.1 shows that,  such entropy solution is the  sharp interface limit for a family global strong solutions of the modified Jin-Xin relaxation system if we choose the thickness of the interface as the relaxation time. In the analysis, it is found that the interaction of shock waves of different families can pass through the interface and maintain the wave strength and wave speed without being affected by the interface for immiscible two-phase flow.
\end{remark}

Here we make some comments on the key steps in the proof. Since the problem considered here involves the perturbations in the vicinity of the shock wave, therefore, it is necessary to use the anti-derivative method to overcome the difficulty that the derivative of shock wave is negative. However, the Navier-Stokes/Allen-Cahn also involves the change of phase field,  the phase field equation does not satisfy the conservation law, the anti-derivative method can't be applied. To resolve this conflict, we adopt the method of using the antiderivative only for the mass equation and momentum equation, while keeping the phase field equation unchanged, see \eqref{Antiderivative equation}. Moreover, in order to get the estimate $\int_0^\tau\|\mathbf{w}\|^2ds$(see \eqref{Energy estimates before the collision}), the weighted estimation method is used(see \eqref{beta}).
Again, in order to apply the anti-derivative method after interaction time,  the relaxation viscous shock profiles with shifts for the outgoing shocks should be constructed(see \eqref{shock after t0}), and the better estimate on the antiderivative of $\widetilde{\mathbf{U}}-\overline{\mathbf{U}}$ is obtained at the interaction time(see \eqref{key1}).

The outline of the proof for Theorem 1.1 is as follows. A priori estimates \eqref{a Priori Estimate Before the Interaction Time} and \eqref{a Priori Estimate After the Interaction Time} of the shock waves of different families and the phase field before and after the interaction time are obtained in section 2,section 3 respectively, then by using $\pmb{\phi}=\mathbf{R}\mathbf{w}$ before the interaction time $t_0$,  and $\tilde{\pmb{\phi}}=\mathbf{\tilde{R}}\mathbf{\tilde{w}}$ after the interaction time $t_0$, so the sharp interface limit is achieved immediately
$\sup_{y\in\mathbb{R}}|\pmb{\phi}_y,\tilde{\pmb{\phi}}_y|\rightarrow0$ as $\epsilon\rightarrow0$. In section 4, the modified Jin-Xin scheme is used to give numerical results for the sharp interface limit of the Navier-Stokes/Allen-Cahn system \eqref{original NSAC} in 1-D.

\section{\normalsize Estimates Before the Interacting Time}
\setcounter{equation}{0}

\ \ \ \ In this section, we consider the stability before the interaction of two reverse shock waves.  Letting
\begin{equation}\label{approximate solution}
\mathbf{U}^{S_1S_2}(y,\tau)=\mathbf{U}^{S_1}(y-s_1\tau)+\mathbf{U}^{S_2}(y-s_2\tau)-\mathbf{U}^*.
\end{equation}
From \eqref{1-relaxation shock wave} and \eqref{2-relaxation shock wave},  $\mathbf{U}^{S_1S_2}(y,\tau)$ satisfying
\begin{eqnarray} \label{incoming shock waves}
\mathbf{U}^{S_1S_2}_{\tau}+\mathbf{F}_y(\mathbf{U}^{S_1S_2})-a^2\mathbf{U}^{S_1S_2}_{yy}+\mathbf{U}^{S_1S_2}_{\tau\tau}=\mathbf{G}_y,
\end{eqnarray}
where
\begin{equation}\label{G}
\mathbf{G}=\mathbf{F}(\mathbf{U}^{S_1S_2})-\mathbf{F}(\mathbf{U}^{S_1})-\mathbf{F}(\mathbf{U}^{S_2})+\mathbf{F}(\mathbf{U}^*).
\end{equation}
Under the scaling transformation \eqref{Scaling transformation},  taking the \eqref{approximate solution}  at $\tau=-\frac{t_0}{\epsilon}$ as the initial value of system (\ref{Relaxation system}), that is
\begin{equation}\label{relaxation system before t0}
\left\{\begin{array}{llll}
\displaystyle \mathbf{U}_\tau+\mathbf{F}_y(\mathbf{U})=(a^2 \mathbf{U}_{yy}-\mathbf{U}_{\tau\tau})+\mathbf{H}_{\tau y}+\mathbf{H}_y, \\
\displaystyle \chi_\tau=-v\mu,\ \ \ \displaystyle \mu=(\chi^3-\chi)-\big(\frac{\chi_y}{v}\big)_y,\\
\displaystyle (\mathbf{U},\chi)\Big|_{\tau=-\frac{t_0}{\epsilon}}=\big(\mathbf{U}^{S_1S_2}(y,-\frac{t_0}{\epsilon}),\chi_0\big).
\end{array}\right.
\end{equation}
Setting
 \begin{equation}\label{Phi}
   \mathbf{\Phi}(y,\tau)=(v-V^{S_1S_2},u-U^{S_1S_2})=\mathbf{U}-\mathbf{U}^{S_1S_2},
 \end{equation}
subtracting the \eqref{incoming shock waves} from \eqref{relaxation system before t0}, one obtains
\begin{eqnarray}
\left\{ \begin{aligned}
&\mathbf{\Phi}_{\tau}+\big(\mathbf{F}(\mathbf{U}^{S_1S_2}+\mathbf{\Phi})-\mathbf{F}(\mathbf{U}^{S_1S_2})\big)_y=a^2\mathbf{\Phi}_{yy}-\mathbf{\Phi}_{\tau\tau}-\mathbf{G}_y+\mathbf{H}_{\tau y}+\mathbf{H}_y,\\
& \chi_\tau=-v\mu,\ \ \ \displaystyle \mu=(\chi^3-\chi)-\Big(\frac{\chi_y}{v}\Big)_y,\\
&(\mathbf{\Phi},\chi)\Big|_{\tau=-\frac{t_0}{\epsilon}}=(\mathbf{0},\chi_0).
\end{aligned} \right.
\end{eqnarray}
The anti-derivative technique should be used here, (cf. \cite{KM85}, \cite{MN85} and the references therein), that is
\begin{equation}\label{anti-derivative technique}
  \pmb{\phi}(y,\tau)=(\phi_1,\phi_2)^\mathbb{T}=\int^y_{-\infty}\mathbf{\Phi}(z,\tau)dz,
\end{equation}
one has
\begin{eqnarray} \label{Antiderivative equation}
\left\{ \begin{aligned}
&\pmb{\phi}_{\tau}+\mathbf{F}'(\mathbf{U}^{S_1S_2})\pmb{\phi}_{y}-a^2\pmb{\phi}_{yy}=-\pmb{\phi}_{\tau\tau}-(\mathbf{G}+\mathbf{Q})+\mathbf{H}_{\tau}+\mathbf{H},\\
& \chi_\tau=-v\mu,\ \ \ \displaystyle \mu=(\chi^3-\chi)-\big(\frac{\chi_y}{v}\big)_y,\\
&\pmb{\phi}\Big|_{\tau=-\frac{t_0}{\epsilon}}=0,
\end{aligned} \right.
\end{eqnarray}
where
\begin{equation}\label{Q}
  \mathbf{Q}=\mathbf{F}(\mathbf{U}^{S_1S_2}+\pmb{\phi}_y)
-\mathbf{F}(\mathbf{U}^{S_1S_2})-\mathbf{F}'(\mathbf{U}^{S_1S_2})\pmb{\phi}_y=O(1)(\pmb{\phi}_y^2),
\end{equation}
and
\begin{eqnarray} \label{sigma}
|\mathbf{G}|\leq C\delta^2e^{-c\delta(|y|+|\tau|)}\overset{\text{def}}{=} \sigma.
\end{eqnarray}
Digitalizing the system (\ref{Antiderivative equation}) by using the following linear transformation
\begin{eqnarray} \label{linear transformation}
\pmb{\phi}=\mathbf{R}(\mathbf{U}^{S_1S_2})\mathbf{w}=(r_{11}w_1+r_{12}w_2,r_{21}w_1+r_{22}w_2),
\end{eqnarray}
and then multiplying the result  by $\mathbf{L}(\mathbf{U}^{S_1S_2})$, one obtains
\begin{eqnarray} \label{Diagonalized system of equations}
\left\{ \begin{aligned}
&\mathbf{w}_{\tau}+\mathbf{\Lambda}\mathbf{w}_y-a^2\mathbf{w}_{yy}=-\mathbf{w}_{\tau\tau}-\mathbf{L}\big(\mathbf{R}_{\tau}-a^2\mathbf{R}_{yy}+\mathbf{F}'\mathbf{R}_y+\mathbf{R}_{\tau\tau}\big)\mathbf{w}\\
&\displaystyle\qquad\qquad\qquad\qquad\qquad\qquad+\mathbf{L}(2a^2\mathbf{R}_y\mathbf{w}_y-\mathbf{G}-\mathbf{Q}+\mathbf{H}_{\tau}+\mathbf{H}-2\mathbf{R}_\tau\mathbf{w}_\tau),\\
& \chi_\tau=-v\mu,\ \ \ \displaystyle \mu=(\chi^3-\chi)-\big(\frac{\chi_y}{v}\big)_y,\\
&(\mathbf{w},\chi)\Big|_{\tau=-\frac{t_0}{\epsilon}}=(\mathbf{0},\chi_0),
\end{aligned} \right.
\end{eqnarray}
where
\begin{equation}\label{L,R,Lambda,v}\begin{split}
 & \mathbf{\Lambda}=\mathbf{\Lambda}(\mathbf{U}^{S_1S_2}),\ \ \mathbf{R}=\mathbf{R}(\mathbf{U}^{S_1S_2}),\ \ \mathbf{L}=\mathbf{L}(\mathbf{U}^{S_1S_2}),\\
 & \mathbf{F}'=\mathbf{F}'(\mathbf{U}^{S_1S_2}),\qquad v=(r_{11}w_1+r_{12}w_2)_y+V^{S_1S_2}.
\end{split}\end{equation}
 Setting the solution space for \eqref{Diagonalized system of equations} as follows:
\begin{equation}\label{Solution space}
\begin{split}
& X_M(I)=\Big\{ (\mathbf{w},\chi)\Big|\mathbf{w}\in C(I;H^3), \mathbf{w}_\tau\in C(I;H^2), \chi^2-1\in C(I;L^2),\chi_y\in C(I;H^3),\\
&\qquad\qquad\qquad\qquad \chi_{y\tau}\in C(I;H^2),\mathbf{w}\in L^2(I;H^3),\mathbf{w}_\tau\in L^2(I;H^2),\\
&\qquad\qquad\qquad\qquad\chi^3-\chi\in L^2(I;L^2),\chi_y\in L^2(I;H^3),\chi^2-1\in L^2(I;L^2),\\
&\qquad\qquad\qquad\qquad\sup_{\tau\in I}\{\|\mathbf{w}\|_3,\|\mathbf{w}_\tau\|_2,\|\chi^2-1\|,\|\chi_y\|_3,\|\chi_{\tau y}\|_1\}\leq M \Big\},
\end{split}\end{equation}
where $I\subseteq [-\frac{t_0}{\epsilon},0]$ is the any  time interval. The proof of the local existence and uniqueness to the system \eqref{Diagonalized system of equations} is trivial,  the detailed omitted.
 From the definition \eqref{Solution space}, choosing $M$ small enough, called as $\delta_0$,  by using  Sobolev embedding theorem, there exist $m_0>0$, such that 
\begin{equation}\label{upper and lower bound for v and chi}
\left.\begin{array}{llll}
\displaystyle0<\frac{3}{4}v_{-}\leq v\leq \frac{5}{4}v_-,\ \inf_{x\in \mathbb{R}, t\in (0,T)}3\chi^2-1\geq m_0>0.
\end{array}\right.\end{equation}
By using the maximum principle for parabolic equation (see Lemma 2.1 in \cite{p2005}) and \eqref{upper and lower bound for v and chi}, we obtain
which yields
\begin{equation}\label{bound of phase field}
|\chi(x,t)|\leq 1,\qquad \forall (x,t)\in (-\infty,+\infty)\times[0,T],
\end{equation}

In the following, the a priori estimates of the solutions for \eqref{Diagonalized system of equations} before the interaction time will be given.

\begin{proposition} \label{a Priori Estimate Before the Interaction Time}
(A Priori Estimate Before the Interaction Time)   Assume that $\mathbf{w}\in X_M[-\frac{t_0}{\epsilon},s]$ is the solution of the system \eqref{Diagonalized system of equations} with $s\leq0$,
 there exist positive constants $\delta_0$, $M_0$ and $C$ independent of $s$, such that, if $\delta\leq\delta_0$ and $M\leq M_0$,
then for any $\tau\in [-\frac{t_0}{\epsilon},s]$, it holds that
\begin{equation}\label{Energy estimates before the collision}
\begin{split}
&\|\mathbf{w}(s)\|_{H^3}^2+\|\mathbf{w}_\tau(s)\|_{H^2}^2+\|\chi^2(s)-1\|^2+\|\chi_{y}(s)\|_3^2+\|\chi_{y\tau}(s)\|_2^2\\
&+\int_{-\frac{t_0}{\epsilon}}^{\tau}\big(\|\mathbf{w}_y(s)\|_{H^2}^2+\|\mathbf{w}_\tau(s)\|_{H^2}^2+\|(\chi^3-\chi)(s)\|^2+\|\chi_y(s)\|_3^2\big)ds\\
&\leq CM_0\Big(\delta^{\frac{1}{2}}e^{\frac{-c\delta|t-t_0|}{\epsilon}}+1\Big),
\end{split}
\end{equation}where $C$ is the positive constants which may depend on $(v_-,u_-)$ but is independent $s$.
\end{proposition}

Proposition \ref{a Priori Estimate Before the Interaction Time} follows from the following lemmas. Lemma 1.1 is the  low-order estimate of the solutions  in the following form:
\begin{lemma} \label{lem of v}
Under the assumption of Proposition \ref{a Priori Estimate Before the Interaction Time}, there exist positive constants $\delta_0$, $M_0$ independent of $T$, such that, if $\delta\leq\delta_0$ and $M\leq M_0$,
then for any $\tau\in [-\frac{t_0}{\epsilon},T]$, it holds that
\begin{equation}\label{the lower estimate of w}
\begin{split}
&\|\mathbf{w}(\tau)\|^2+\|\mathbf{w}_{\tau}(\tau)\|^2+\|\chi^2-1\|^2+\int_{-\frac{t_0}{\epsilon}}^{\tau}\int_\mathbb{R}\big(|\mathbf{U}^{S_1}_y|+|\mathbf{U}^{S_2}_y|\big)\big(w_1^2+w_2^2\big)dyds\\
&\quad+\int_{-\frac{t_0}{\epsilon}}^{\tau}\big(\|\mathbf{w}_y(s)\|^2+\|\mathbf{w}_{\tau}(s)\|^2+\|\chi^3-\chi\|^2+\|\chi_y\|^2\big)ds\leq CM_0\big(\delta^{\frac{1}{2}}e^{-c\delta|\tau|}+1\big),
\end{split}
\end{equation}where $C$ is the positive constants which may depend on $(v_-,u_-)$ but is independent $T$.
\end{lemma}

\begin{proof} By using the maximum principle for parabolic equation (see Lemma 2.1 in \cite{p2005}) and \eqref{upper and lower bound for v and chi}, we obtain
\begin{equation}\label{phi333}
\displaystyle \chi^2\leq1.
\end{equation}
Rewriting \eqref{Diagonalized system of equations}$_2$  as
\begin{equation}\label{rewrite of phase field equation}
 \chi_\tau=-v(\chi^3-\chi)+v\big(\frac{\chi_y}{v}\big)_y.
\end{equation}
In order to get the estimate $\int_0^\tau\|\mathbf{w}\|^2ds$, thanks for the weighted function introduced by Xin \cite{Xi00}, the weighted estimation method is used below. Setting
\begin{equation}\label{alpha}
\alpha_1=1+\beta_1^2,\qquad\alpha_2=1+\beta_2^1,
\end{equation}
where, $\beta_1^2$ and  $\beta_2^1$ are expressed as the following expression
\begin{equation}\label{beta}
\begin{split}
\beta_1^2(\mathbf{U}^{S_2})=\frac{\lambda_1(\mathbf{U}^{S_2}(0))-s_2}{\lambda_1(\mathbf{U}^{S_2})-s_2}e^{-\delta^{-\frac14}\int_0^{\xi_2}\frac{|\lambda_{2y}(\mathbf{U}^{S_2})|}{\lambda_1(\mathbf{U}^{S_2})-s_2}dy},\\
\beta_2^1(\mathbf{U}^{S_1})=\frac{\lambda_2(\mathbf{U}^{S_1}(0))-s_1}{\lambda_2(\mathbf{U}^{S_1})-s_1}e^{-\delta^{-\frac14}\int_0^{\xi_1}\frac{|\lambda_{1y}(\mathbf{U}^{S_1})|}{\lambda_2(\mathbf{U}^{S_1})-s_1}dy}.
\end{split}
\end{equation}
By direct calculation, one has
\begin{equation}
\partial_{\xi_j}\big\{[\lambda_i(\mathbf{U}^{S_j})-s_j]\beta_i^j(\mathbf{U}^{S_j})\big\}=-\delta^{-\frac14}|\partial_{\xi_j}\lambda_j(\mathbf{U}^{S_j})|\beta_i^j(\mathbf{U}^{S_j}),\ j\neq i,\ i,j=1,2,
\end{equation}
moreover, according to the eigenvalue properties of the strict hyperbolic system, there exists a positive constant $C$, such that
 \begin{equation}\label{Upper and lower bounds for alpha}
  C^{-1}<\alpha_i<C.
 \end{equation}
Setting
$$\boldsymbol{\alpha}=\mathrm{diag}(\alpha_1,\alpha_2),$$
multiplying (\ref{Diagonalized system of equations})$_1$ by $\mathbf{w}^\mathbb{T}\pmb{\alpha}$, \eqref{rewrite of phase field equation} by $\chi^3-\chi$, adding them together,  integrating the result with respect to $y$ over $\mathbb{R}$, yields
\begin{equation}
\begin{split}
&\frac12\frac d{d\tau}\int_{\mathbb{R}}\big(\mathbf{w}^{\mathbb{T}}\boldsymbol{\alpha} \mathbf{w}+\frac{(\chi^2-1)^2}{2}\big)dy\\
&\qquad+\int_{\mathbb{R}}\big( a^2\mathbf{w}^{\mathbb{T}}_y\boldsymbol{\alpha}\mathbf{w}_{y}+v(\chi^3-\chi)^2+(3\chi^2-1)\chi_y^2\big)dy\\
&=I_1+I_2+I_3+I_4+I_5,
\end{split}
\end{equation}
and
\begin{eqnarray}
&&I_1=\frac{1}{2}\int_{\mathbb{R}}\sum_{i=1}^2\Big(\alpha_{i\tau}+\big[\alpha_i\lambda_i(\mathbf{U}^{S_1S_2})\big]_y+\big(a^2- s_i^2\big)\alpha_{iyy}\Big) w_i^2dy,\\
&&I_2=\int_{\mathbb{R}}\mathbf{w}^\mathbb{T}\boldsymbol{\alpha} \mathbf{L}(2a^2\mathbf{R}_y\mathbf{w}_y-\mathbf{G}-\mathbf{Q}-2\mathbf{R}_\tau \mathbf{w}_\tau)dy,\\
&&I_3=-\int_{\mathbb{R}}\mathbf{w}^\mathbb{T}\boldsymbol{\alpha} \mathbf{L}(\mathbf{R}_{\tau}-a^2\mathbf{R}_{yy}+\mathbf{F}'(\mathbf{U}^{S_1S_2})\mathbf{R}_y+\mathbf{R}_{\tau\tau})\mathbf{w}dy,\\
&&I_4=-\int_{\mathbb{R}}\frac{\big((r_{11}w_1+r_{12}w_2)_y+V^{S_1S_2}\big)_y}{v}\chi_y(\chi^3-\chi)dy,\\
&&I_5=\int_{\mathbb{R}}\mathbf{w}^\mathbb{T}\boldsymbol{\alpha} \mathbf{L}(\mathbf{H}_{\tau}+\mathbf{H})dy.
\end{eqnarray}
In the following, we give the estimates of $I_1,I_2,I_3,I_4,I_5$ in turn. To do this, we start by estimating the derivatives of the eigenvalues:
\begin{eqnarray}
\lambda_{1y}(\mathbf{U}^{S_1S_2})&=&\lambda'_1(\mathbf{U}^{S_1S_2})(\mathbf{U}^{S_1}_y+\mathbf{U}^{S_2}_y) \notag\\
&=&\lambda_{1y}(\mathbf{U}^{S_1})+(\lambda'_1(\mathbf{U}^{S_1S_2})
-\lambda'_1(\mathbf{U}^{S_1}))\mathbf{U}^{S_1}_y+\lambda'_1(\mathbf{U}^{S_1S_2})\mathbf{U}^{S_2}_y\notag\\
&=&\lambda_{1y}(\mathbf{U}^{S_1})+O(1)(\sigma+|\mathbf{U}^{S_2}_y|),
\end{eqnarray}
where $\sigma$  is the positive constant given in \eqref{sigma},  similarly
\begin{eqnarray}
\lambda_{2y}(\mathbf{U}^{S_1S_2})
=\lambda_{2y}(\mathbf{U}^{S_2})+O(1)(\sigma+|\mathbf{U}^{S_1}_y|),
\end{eqnarray}
then\begin{eqnarray}
&&\alpha_{1\tau}+\big(\alpha_1\lambda_1(\mathbf{U}^{S_1S_2})\big)_y=\alpha_{1\tau}+\lambda_1(\mathbf{U}^{S_1S_2})\alpha_{1y}+\alpha_1\lambda_{1y}(\mathbf{U}^{S_1S_2})\notag\\
&&=(\lambda_1(\mathbf{U}^{S_1S_2})-s_2)\beta_{1y}^2(\mathbf{U}^{S_2})+\alpha_1\lambda_{1y}(\mathbf{U}^{S_1S_2})\notag\\
&& =\big(\lambda_1(\mathbf{U}^{S_2})-s_2\big)\beta_{1y}^2(\mathbf{U}^{S_2})+\big(\lambda_1(\mathbf{U}^{S_1S_2})-\lambda_1(\mathbf{U}^{S_2})\big)\beta_{1y}^2(\mathbf{U}^{S_2})+\alpha_1\lambda_{1y}(\mathbf{U}^{S_1S_2})\notag\\
&&=-O(1)(\delta^{-\frac14}|\mathbf{U}^{S_2}_y|+\sigma)+\alpha_1\lambda_{1y}(\mathbf{U}^{S_1})+O(1)\sigma,
\end{eqnarray}
and by the same way, one has
\begin{eqnarray}
&&\alpha_{2\tau}+\big(\alpha_2\lambda_2(\mathbf{U}^{S_1S_2})\big)_y=-O(1)(\delta^{-\frac14}|\mathbf{U}^{S_1}_y|+\sigma)+\alpha_2\lambda_{2y}(\mathbf{U}^{S_2})+O(1)\sigma.
\end{eqnarray}
Moreover, one obtains
\begin{eqnarray}
&&\int_\mathbb{R}\sum_{i=1}^{2}(a^2- s_i^2)\alpha_{iyy}w_i^2dy \leq O(1)\delta^{\frac12}\int_\mathbb{R}(|\mathbf{U}^{S_2}_y|w_1^2+|\mathbf{U}^{S_1}_y|w_2^2)dy.
\end{eqnarray}
So, one achieves
\begin{eqnarray}\label{estimate of J}
&&I_1\leq - O(1)\delta^{-\frac14}\int_{\mathbb{R}}(|\mathbf{U}^{S_2}_y|w_1^2 + |\mathbf{U}^{S_1}_y|w_2^2)dy+\int_{\mathbb{R}}\big(\alpha_1|\lambda_{1y}(\mathbf{U}^{S_1})|w_1^2+\alpha_2|\lambda_{2y}(\mathbf{U}^{S_2})|w_2^2\big)dy\notag\\
&&\qquad\qquad+O(1)\int_{\mathbb{R}}\sigma\mathbf{w}^2dy.
\end{eqnarray}
Further, the estimate of $I_2$ is as follows
\begin{eqnarray}\label{estimate of I2}
I_2&\leq&O(1)\int_{\mathbb{R}}\big[|\mathbf{w}|(|\mathbf{U}^{S_1S_2}_y||\mathbf{w}_y|+|\mathbf{U}^{S_1S_2}_y\mathbf{w}|^2+|\mathbf{w}_y|^2+\sigma)+|\mathbf{w}||\mathbf{U}^{S_1S_2}_y||\mathbf{w}_\tau|\big]dy\notag\\
&\leq& O(1)(M_0+\delta^{\frac{1}{2}})\int_{\mathbb{R}}|\mathbf{w}_y|^2dy+O(1)\delta^{\frac32}\int_{\mathbb{R}}(|\mathbf{U}_y^{S_1}|+|\mathbf{U}_y^{S_2}|)\mathbf{w}^2dy\notag\\
&&+O(1)\int_{\mathbb{R}}\sigma|\mathbf{w}|dy+O(1)\delta^{\frac{1}{2}}\int_{\mathbb{R}}|\mathbf{w}_\tau|^2dy.
\end{eqnarray}
 Now, we consider the $I_3$. Letting
 \begin{equation}\label{a-ij}
  \begin{split}
   &(a_{ij})_{2\times2}=\mathbf{L}(\mathbf{U}^{S_1S_2})\Big(\mathbf{R}_{\tau}(\mathbf{U}^{S_1S_2})-a^2\mathbf{R}_{yy}(\mathbf{U}^{S_1S_2})\\
   &\qquad\qquad\qquad\qquad\qquad+\mathbf{F}'(\mathbf{U}^{S_1S_2})\mathbf{R}_y(\mathbf{U}^{S_1S_2})+\mathbf{R}_{\tau\tau}(\mathbf{U}^{S_1S_2})\Big),
\end{split} \end{equation}
calculating directly, one has
\begin{eqnarray}\label{estimate of I1}
I_3&=&\Big|\int_{\mathbb{R}}\sum_{i,j=1}^2a_{ij}w_iw_jdy\Big|\notag\\
&\leq&\int_{\mathbb{R}}O(1)(1+\eta^{-1})(|\mathbf{U}^{S_2}_y|w_1^2+|\mathbf{U}^{S_1}_y|w_2^2)dy\\
&&+\int_{\mathbb{R}}O(1)(\eta+\delta)(|\mathbf{U}^{S_1}_y|w_1^2+|\mathbf{U}^{S_2}_y|w_2^2)+O(1)\sigma|\mathbf{w}|^2dy,\notag
\end{eqnarray}
where we choose $\eta$ so small such that $O(1)(\eta+\delta)|\mathbf{U}^{S_i}_y|\leq\alpha_i|\lambda(\mathbf{U}^{S_i})_y|$, $i=1,2$ .
About the $I_4$, one has
\begin{equation}\label{I4}
\begin{split}
 I_4&\leq \int_{\mathbb{R}}\big|\big((r_{11}w_1+r_{12}w_2)_{yy}+V_y^{S_1S_2})\big|\chi_y\big|\big|\chi^3-\chi\big|dy\\
 &\leq C\|\mathbf{w}_{yy}\|^{\frac12}\|\mathbf{w}_{yyy}\|^{\frac12}\big(\int_{\mathbb{R}}\chi_y^2dy+\int_{\mathbb{R}}(\chi^3-\chi)^2dy\big)\\
  &\leq CM\big(\int_{\mathbb{R}}\chi_y^2dy+\int_{\mathbb{R}}(\chi^3-\chi)^2dy\big).
 \end{split}
\end{equation}
For $I_5$, since
\begin{eqnarray}
  \Big(\frac{\chi_y}{v}\Big)^2_\tau&=&2 \frac{\chi_y}{v}\Big(\frac{\chi_y}{v}\Big)_\tau \notag\\
  &=&-2\frac{\chi_y}{v^2}v_y(\chi^3-\chi)-2\frac{\chi^2_y}{v}(3\chi^2-1)+2\frac{\chi_y\chi_{yy}v_y}{v^3}-\frac{\chi_y^2v_y^2}{v^3}+ 2\frac{\chi_y\chi_{yyy}}{v^2}\\
   &&-2\frac{\chi_{yy}\chi_yv_y}{v^3}-2\frac{\chi_{yy}v_y\chi_y}{v^3}-2\frac{\chi_{y}^2v_{yy}}{v^3}-2\frac{\chi_y^2v_y^2}{v^3},\notag
\end{eqnarray}
one gets
\begin{equation}\label{I5}
 I_5= \int_{\mathbb{R}}\mathbf{w}^\mathbb{T}\boldsymbol{\alpha} \mathbf{L}(\mathbf{H}_{\tau}+\mathbf{H})dy\leq CM\big(\|\mathbf{w}\|^2+\|\chi_{y}\|^2\big),
\end{equation}
and thus, choosing  $\delta,M$ and $\eta_0$ small enough, combining with (\ref{estimate of J}), (\ref{estimate of I2}) and (\ref{estimate of I1}), \eqref{I4}, \eqref{I5}, one obtains
\begin{eqnarray}\label{estimate 1}
&&\frac d{d\tau}\int_{\mathbb{R}}\big(\mathbf{w}^{\mathbb{T}}\boldsymbol{\alpha} \mathbf{w}+\frac{(\chi^2-1)^2}{4}\big)dy+O(1)\delta^{-\frac14}\int_{\mathbb{R}}(|\mathbf{U}^{S_2}_y|w_1^2+|\mathbf{U}^{S_1}_y|w_2^2)dy
\notag\\
&&+O(1)\int_{\mathbb{R}}\big( a^2\mathbf{w}^{\mathbb{T}}_y\boldsymbol{\alpha}\mathbf{w}_{y}+v(\chi^3-\chi)^2+(3\chi^2-1)\chi_y^2\big)dy\\
&&\leq O(1)\int_{\mathbb{R}}\sigma|\mathbf{w}|dy+O(1)\delta^{\frac{1}{2}}\int_{\mathbb{R}}|\mathbf{w}_\tau|^2dy.\notag
\end{eqnarray}
Finally, multiply (\ref{Diagonalized system of equations})$_1$ by $\mathbf{w}^\mathbb{T}_{\tau}$, integrate with respect to $y$ over $\mathbb{R}$, we have
\begin{equation}
 \frac{1}{2}\frac{ d}{d\tau}\int_\mathbb{R}\big(a^2\mathbf{w}^{\mathbb{T}}_y\mathbf{w}_{y}+\mathbf{w}^{\mathbb{T}}_\tau\mathbf{w}_{\tau}\big)dy+\int_\mathbb{R}\big(\mathbf{w}^{\mathbb{T}}_\tau\mathbf{w}_{\tau}+\mathbf{w}^{\mathbb{T}}_\tau\mathbf{\Lambda}\mathbf{w}_{y}
\big)dy=I_5+I_6+I_7,
\end{equation}
where
\begin{eqnarray}
&&I_6=-\int_{\mathbb{R}}\mathbf{w}^\mathbb{T}_{\tau} \mathbf{L}\big(\mathbf{R}_{\tau}-a^2\mathbf{R}_{yy}+\mathbf{F}'(\mathbf{U}^{S_1S_2})\mathbf{R}_y+\mathbf{R}_{\tau\tau}\big)\mathbf{w}dy,\\
&&I_7=\int_{\mathbb{R}}\mathbf{w}_{\tau}^\mathbb{T} \mathbf{L}(2a^2\mathbf{R}_y\mathbf{w}_y-\mathbf{G}-\mathbf{Q}-2\mathbf{R}_\tau \mathbf{w}_\tau)dy,\\
&&I_8=\int_{\mathbb{R}}\mathbf{w}_\tau^\mathbb{T} \mathbf{L}(\mathbf{H}_{\tau}+\mathbf{H})dy.
\end{eqnarray}
Then one gets
\begin{equation}
\begin{split}
&I_6\leq O(1)\delta^{\frac{1}{2}}\int_{\mathbb{R}}\!\big(|\mathbf{U}^{S_1}_y|+|\mathbf{U}^{S_2}_y|\big)\big(w_1^2
+w_2^2\big)dy+O(1)\big(\delta^{\frac12}\|\mathbf{w}_{\tau}\|^2+\int_{\mathbb{R}}\sigma|\mathbf{w}|dy\big),
\end{split}
\end{equation}
and
\begin{equation}\begin{split}
&I_7\leq O(1)\delta^{\frac{1}{2}}\int_{\mathbb{R}}\!\big(|\mathbf{U}^{S_1}_y|+|\mathbf{U}^{S_2}_y|\big)\big(w_1^2
+w_2^2\big)dy\\
&\qquad\quad+O(1)(\delta^{\frac{1}{2}}+\eta_0)\big(\|\mathbf{w}_y\|^2+\|\mathbf{w}_{\tau}\|^2\big)+O(1)\int_\mathbb{R}\sigma|\mathbf{w}_{\tau}|dy.
\end{split}\end{equation}
Moreover, for $I_8$, one has
\begin{equation}\begin{split}
&I_8= \int_{\mathbb{R}}\mathbf{w}^\mathbb{T}_\tau\boldsymbol{\alpha} \mathbf{L}(\mathbf{H}_{\tau}+\mathbf{H})dy\leq CM\big(\|\mathbf{w}_\tau\|^2+\|\chi_{y}\|^2\big),
\end{split}\end{equation}
thus, it holds that
\begin{equation}\label{estimate 2}
\begin{split}
  & \frac{1}{2}\frac{ d}{d\tau}\int_\mathbb{R}\big(a^2\mathbf{w}^{\mathbb{T}}_y\mathbf{w}_{y}+\mathbf{w}^{\mathbb{T}}_\tau\mathbf{w}_{\tau}\big)dy+\int_\mathbb{R}\big(\mathbf{w}^{\mathbb{T}}_\tau\mathbf{w}_{\tau}+\mathbf{w}^{\mathbb{T}}_\tau\mathbf{\Lambda}\mathbf{w}_{y}\big)dy\\
& \leq O(1)\delta^{\frac{1}{2}}\int_{\mathbb{R}}\!\big(|\mathbf{U}^{S_1}_y|+|\mathbf{U}^{S_2}_y|\big)\big(w_1^2
+w_2^2\big)dy+O(1)\int_\mathbb{R}\sigma\big(|\mathbf{w}|+|\mathbf{w}_\tau|\big)dy\\
&\qquad+O(1)(\delta^{\frac12}+\eta_0)\big(\|\mathbf{w}_y\|^2+\|\mathbf{w}_\tau\|^2\big)+CM\big(\|\mathbf{w}_\tau\|^2+\|\chi_{y}\|^2\big).
\end{split}\end{equation}
Therefore, adding (\ref{estimate 1}) and (\ref{estimate 2}), integrating \eqref{estimate 2} with respect to $s$ over $[-\frac{\epsilon}{t_0}, \tau]$ , by using the sub-characteristic condition (\ref{sub-characteristic condition}) and the definition of $\alpha_i$,
we  obtain the following lower order estimate
\begin{equation}\label{estimate of v}
\begin{split}
&\|\mathbf{w}(\tau)\|^2+\|\mathbf{w}_{\tau}(\tau)\|^2+\int_{-\frac{t_0}{\epsilon}}^{\tau}\int_\mathbb{R}\big(|\mathbf{U}^{S_1}_y|+|\mathbf{U}^{S_2}_y|\big)\big(w_1^2+w_2^2\big)dyds\\
&+\int_{-\frac{t_0}{\epsilon}}^{\tau}\big(\|\mathbf{w}_y(s)\|^2+\|\mathbf{w}_{\tau}(s)\|^2\big)ds
\leq O(1)\Big(\int_{-\frac{t_0}{\epsilon}}^{\tau}\int_\mathbb{R}\sigma\big(|\mathbf{w}|+|\mathbf{w}_\tau|\big)dy+\|\chi_0^2-1\|^2\Big).\end{split}
\end{equation}
In the following, we will give the estimate of the right-hand side of the inequality (\ref{estimate of v}) as follows
\begin{eqnarray}
&&\int_{-\frac{t_0}{\epsilon}}^{\tau}\int_{\mathbb{R}}\sigma|\mathbf{w}|dyds=O(1)\delta^2\int_{-\frac{t_0}{\epsilon}}^{\tau}\int_{\mathbb{R}}e^{-c\delta|y|-
c\delta|s|}|\mathbf{w}|dyds\notag\\
&&\leq O(1)\delta^{\frac{3}{2}}\int_{-\frac{t_0}{\epsilon}}^{\tau}e^{-c\delta|s|}\|\mathbf{w}\|ds\leq O(1)M_0\delta^{\frac{1}{2}}e^{-c\delta|\tau|},
\end{eqnarray}
noticing $|\mathbf{w}_{\tau}|\leq O(1)|\mathbf{w}_{y}|$, by the same way, we get
\begin{equation}
\int_{-\frac{t_0}{\epsilon}}^{\tau}\int_{\mathbb{R}}\sigma|\mathbf{w}_\tau|dyds\leq  O(1)M_0\delta^{\frac{1}{2}}e^{-c\delta|\tau|}.
\end{equation}
Thus, the   proof of the Lemma \ref{lem of v} is completed.
\end{proof}

\

In the following Lemma, the higher-order estimates of the solutions will be given.
\begin{lemma} \label{lem of high order estimate}
Under the assumption of Proposition \ref{a Priori Estimate Before the Interaction Time}, there exist positive constants $\delta_0$, $M_0$ independent of $T$, $T\leq0$, such that, if $\delta\leq\delta_0$ and $M\leq M_0$,
then for any $s\in [-\frac{t_0}{\epsilon},T]$, it holds that
\begin{equation}\label{Higher Derivatives}
\begin{split}
&\|\mathbf{w}_{y}(s)\|_2^2+\|\mathbf{w}_{\tau  y}(s)\|_1^2+\|\chi_{y}(s)\|_3^2+\|\chi_{\tau y}(s)\|_2^2\\
&+\int_{\frac{-t_0}{\epsilon}}^{s}\Big(\|\mathbf{w}_{y}(s)\|_2^2+\|\mathbf{w}_{y\tau}(s)\|_1^2+\|\chi_{y}(s)\|_2^2\Big)ds\leq C\Big(\delta^{\frac{1}{2}}e^{-\frac{c\delta|t-t_0|}{\epsilon}}+M_0\Big).
\end{split}
\end{equation}
\end{lemma}
\begin{proof}
In order to get the higher order estimate of the solution for the system (\ref{Diagonalized system of equations}), applying $\partial_y$ to system (\ref{Diagonalized system of equations}),  one has
\begin{equation}\label{Higher-order system}
\left\{
\begin{split}
&\mathbf{w}_{\tau y}+(\mathbf{\Lambda}\mathbf{w}_y)_y-a^2\mathbf{w}_{yyy}+\mathbf{w}_{\tau\tau y}=-\big\{\mathbf{L}\big(\mathbf{R}_{\tau}-a^2\mathbf{R}_{yy}+\mathbf{F}'\mathbf{R}_y+\mathbf{R}_{\tau\tau}\big)\mathbf{w}\big\}_y\\
&\qquad\qquad\qquad\qquad\qquad+\big\{\mathbf{L}(2a^2\mathbf{R}_y\mathbf{w}_y-\mathbf{G}-\mathbf{Q}+\mathbf{H}_{\tau}+\mathbf{H}-2\mathbf{R}_\tau\mathbf{w}_\tau)\big\}_y,\\
&\chi_{\tau y}=-\big(v(\chi^3-\chi)\big)_y+\Big(v\big(\frac{\chi_y}{v}\big)_y\Big)_y,\\
&(\mathbf{w},\chi)\Big|_{\tau=-\frac{t_0}{\epsilon}}=(\mathbf{0},\chi_0).
\end{split}\right.
\end{equation}
Multiplying \eqref{Higher-order system}$_1$ by $\eta\mathbf{w}_y^\mathbb{T}$,  \eqref{Higher-order system}$_2$ by $\chi_y$, integrating the resultant over $\mathbb{R}$, yields that
\begin{equation}\label{Higher order derivative inequality 1}
\begin{split}
&\frac12\frac d{d\tau}\int_{\mathbb{R}}\big(\eta|\mathbf{w}_y|^2+2\eta\mathbf{w}_y^\mathbb{T}\mathbf{w}_{y\tau}+\chi_y^2\big)dy+\int_{\mathbb{R}}\big(a^2\eta|\mathbf{w}_{yy}|^2+v(3\chi^2-1)\chi_y^2+\chi^2_{yy}\big)dy\\
& =\eta\int_{\mathbb{R}}\mathbf{w}_{yy}^\mathbb{T}\big[\mathbf{\Lambda}\mathbf{w}_y+\mathbf{L}\big(\mathbf{R}_{\tau}-a^2\mathbf{R}_{yy}+\mathbf{F}'\mathbf{R}_y+\mathbf{R}_{\tau\tau}
\mathbf{w}\big)\big]dy\\
&\qquad-\eta\int_{\mathbb{R}}\mathbf{w}_{yy}^\mathbb{T}\mathbf{L}(2a^2\mathbf{R}_y\mathbf{w}_y-\mathbf{G}-\mathbf{Q}+\mathbf{H}_{\tau}+\mathbf{H}-2\mathbf{R}_\tau \mathbf{w}_\tau)dy\\
&\qquad+\eta\int_{\mathbb{R}}|\mathbf{w}_{y\tau}|^2dy+\int_{\mathbb{R}}\frac{v_y\chi_y}{v}\chi_{yy}-v_y(\chi^3-\chi)\chi_ydy,
\end{split}
\end{equation}
where $\eta$ is a positive constant to be determined later. Furthermore, multiplying (\ref{Higher-order system})$_1$ by $\mathbf{w}_{y\tau}^{\mathbb{T}}$, \eqref{Higher-order system}$_2$ by $\chi_{y\tau}$, integrating over $\mathbb{R}$ with respect to $y$, one gets
\begin{equation}\label{Higher order derivative inequality 2}
\begin{split}
&\frac12\frac d{d\tau}\int_{\mathbb{R}}\big(a^2|\mathbf{w}_{yy}|^2+|\mathbf{w}_{y\tau}|^2+v(3\chi^2-1)\chi_y^2+\chi_{yy}^2+\chi_{y\tau}^2\big)dy\\
&\qquad+\int_{\mathbb{R}}\big(\mathbf{w}_{y\tau}^{\mathbb{T}}\mathbf{\Lambda}\mathbf{w}_{yy}+\mathbf{w}_{y\tau}^2+\|\chi_{yy}\|^2\big)dy\\
& =\int_{\mathbb{R}}\mathbf{w}_{y\tau}^\mathbb{T}\big\{-\mathbf{\Lambda}_y\mathbf{w}_y-\mathbf{L}\big(\mathbf{R}_{\tau}-a^2\mathbf{R}_{yy}+\mathbf{F}'
\mathbf{R}_y+\mathbf{R}_{\tau\tau}\mathbf{w}\big)_y\big\}\\
&\qquad+\int_{\mathbb{R}}\mathbf{w}_{y\tau}^\mathbb{T}\mathbf{L}(2a^2\mathbf{R}_y\mathbf{w}_y-\mathbf{G}-\mathbf{Q}+\mathbf{H}_{\tau}+\mathbf{H}-2\mathbf{R}_\tau \mathbf{w}_\tau)_ydy\\
&\qquad+\int_{\mathbb{R}}\Big(\frac12u_y(3\chi^2-1)\chi^2_y+3v^2\big[\big(\frac{\chi_y}{v}\big)_y-(\chi^3-\chi)\big]\chi_y^2+\frac{\chi_yv_y}{v}\chi_{yy\tau}\Big)dy.
\end{split}
\end{equation}
Integrating \eqref{Higher order derivative inequality 1}, \eqref{Higher order derivative inequality 2}  over $[-\frac{t_0}{\epsilon},s]$ respectively, adding the two integrals together,  for $M$, $\delta$, $\eta$ small enough, one obtains
\begin{equation}\label{Higher order derivative estimation}
\begin{split}
&\|\mathbf{w}_{y}(s)\|_1^2+\|\mathbf{w}_{y\tau}(s)\|^2+\|\chi_y(s)\|_1^2+\|\chi_{y\tau}(s)\|^2\\
&+\int_{\frac{-t_0}{\epsilon}}^{\tau}\big(\|\mathbf{w}_{yy}\|^2+\|\mathbf{w}_{y\tau}\|^2+\|\chi_{yy}\|^2\big)dy \leq C\Big(\delta^{\frac{1}{2}}e^{-\frac{c\delta|t-t_0|}{\epsilon}}+M_0\Big).
\end{split}\end{equation}
In a similar way, we further calculate the higher derivative of the system (\ref{Higher-order system}), and make the corresponding analysis, then the energy estimate of the higher derivative be achieved as following,
\begin{equation}\label{Estimation of higher derivatives}
\begin{split}
&\|\mathbf{w}_{yyy}(s)\|^2+\|\mathbf{w}_{yy\tau}(s)\|^2+\|\chi_{yyy}(s)\|_1^2+\|\chi_{yy\tau}(s)\|_1^2\\
&+\int_{\frac{-t_0}{\epsilon}}^{\tau}\big(\|\mathbf{w}_{yyy}\|^2+\|\mathbf{w}_{yy\tau}\|^2+\|\chi_{yyy}\|_1^2\big)dy \leq C\Big(\delta^{\frac{1}{2}}e^{-\frac{c\delta|t-t_0|}{\epsilon}}+M_0\Big).
\end{split}
\end{equation}
So the proof of the Proposition \ref{a Priori Estimate Before the Interaction Time} is completed.
\end{proof}

\section{\normalsize Estimates After the Interacting Time}
\setcounter{equation}{0}
\ \ \ \ In this section, we discuss the relaxation system for immiscible two-phase flow after the interaction time $\tau=0$ as following
\begin{equation}\label{visous conservation law in section 4}
\left\{\begin{array}{llll}
\displaystyle \mathbf{U}_\tau+\mathbf{F}_y(\mathbf{U})=(a^2 \mathbf{U}_{yy}-\mathbf{U}_{\tau\tau})+\mathbf{H}_{\tau y}+\mathbf{H}_y, \\
\displaystyle \chi_\tau=-v\mu,\ \ \ \displaystyle \mu=(\chi^3-\chi)-\big(\frac{\chi_y}{v}\big)_y,
\displaystyle
\end{array}\right.
\end{equation}
with initial data which has already been solved in section 2.
In order to use the anti-derivative technique, the shock profiles should be modified with suitable shifts. To do this, letting
\begin{equation}\label{shift}
  (\pmb{\nu}_1,\pmb{\nu}_2)=(\mathbf{\tilde{U}}^*-\mathbf{U}_-,\mathbf{U}_+-\mathbf{\tilde{U}}^*).
\end{equation}
It is easy to know that the two vectors $\pmb{\nu}_1,\pmb{\nu}_2$ are linearly independent, so the mass at $\tau=0$
\begin{eqnarray} \mathbf{I}_0 \overset{\text{def}}{=}\int_{\mathbb{R}}
(\mathbf{U}-\mathbf{U}^{\tilde{S}_1\tilde{S}_2})(y,0)dy,
\end{eqnarray}can be represented linearly by $\pmb{\nu}_1,\pmb{\nu}_2$, where $\mathbf{U}^{\tilde{S}_1\tilde{S}_2}(y,\tau)=\mathbf{U}^{\tilde{S}_1}(y-\tilde{s}_1\tau)+\mathbf{U}^{\tilde{S}_2}(y-\tilde{s}_2\tau)-\mathbf{\tilde{U}}^*$. That is,  there exist constants   $b_1$, $b_2$, such that, $\mathbf{I}_0=b_1\pmb{\nu}_1+b_2\pmb{\nu}_2$, furthermore,
we define the relaxation viscous shock profiles with shifts for the outgoing shocks
\begin{eqnarray} \label{shock after t0}
\tilde{\mathbf{U}}^{\tilde{S}_1\tilde{S}_2}(y,\tau)=\mathbf{U}^{\tilde{S}_1}(y-\tilde{s}_1\tau+b_1)+\mathbf{U}^{\tilde{S}_2}(y-\tilde{s}_2\tau+b_2)-\mathbf{\tilde{U}}^*.
\end{eqnarray}
By a simple calculation,  it yields that
\begin{eqnarray}\label{initial data after t0}
\mathbf{I}(b_1,b_2)=\int_{\mathbb{R}}(\mathbf{U}-\tilde{\mathbf{U}}^{\tilde{S}_1\tilde{S}_2})(y,0)dy=\mathbf{I}_0+\frac{\partial \mathbf{I}}{\partial b_1}\cdot
b_1+\frac{\partial \mathbf{I}}{\partial b_2}\cdot b_2=0.
\end{eqnarray}
As in the previous section, we also need the smallness of the antiderivative of perturbation around $\widetilde{\mathbf{U}}$ at $\tau=0$,
  which seems not obvious. Different from \cite{HWWY} and \cite{SZ}, we need the following key Lemma first.
 \begin{lemma} \label{key lemma1}
 For the relaxation  shock waves of system \eqref{visous conservation law in section 4}, the following estimation holds
 \begin{eqnarray}
 &&|v^{S_i}_{yy}|=O(1)\delta_i|v^{S_i}_y|,\,\,\,|u^{S_i}_{yy}|=O(1)\delta_i|u^{S_i}_y|, \ \ i=1,2,\notag\\
 &&|v^{\widetilde{S}_i}_{yy}|=O(1)\tilde{\delta}_i|v^{\widetilde{S}_i}_y|,\,\,\,|u^{\widetilde{S}_i}_{yy}|
 =O(1)\tilde{\delta}_i|u^{\widetilde{S}_i}_y|,\ \  i=1,2.\notag
\end{eqnarray}
 \end{lemma}
\begin{proof}Thanks for the results in \cite{Xi00}, that $|\mathbf{U}^{S_1}_{yy}|=O(1)\delta_1|\mathbf{U}^{S_1}_{y}|$ holds, so it is only need to prove $v^{S_1}_y\sim u^{S_1}_y$. By using (\ref{1-relaxation shock wave}), one has
\begin{eqnarray}
\left\{ \begin{aligned}
 &-s_1v^{S_1}_y-u^{S_1}_y\sim O(1)\delta_1|v^{S_1}_y|,\notag\\
 &-s_1u^{S_1}_y+p_vv^{S_1}_y\sim O(1)\delta_1|u^{S_1}_y|,\notag
\end{aligned} \right.
\end{eqnarray}
and  $|v^{S_1}_{yy}|=O(1)\delta_1|v^{S_1}_y|$, $|u^{S_1}_{yy}|=O(1)\delta_1|u^{S_1}_y|$ are immediately available. Other relations in this Lemma can be proved similarly, here we omit.
\end{proof}

The following is the difference between the relaxation viscous shock waves before and after the interaction time,
\begin{eqnarray} &&(\tilde{\mathbf{U}}^{\tilde{S}_1\tilde{S}_2}-\mathbf{U}^{S_1S_2})\Big|_{\tau=0}=
\left\{ \begin{aligned}
&[(\mathbf{U}^{\tilde{S}_1}(y+b_1)-\mathbf{U}_-)-(\mathbf{U}^{S_1}(y)-\mathbf{U}^*)]\notag\\
&+[(\mathbf{U}^{\tilde{S}_2}(y+b_2)-\mathbf{\tilde{U}}^*)-(\mathbf{U}^{S_2}(y)-\mathbf{U}_-)], \,\,\,y\leq0,\\
&[(\mathbf{U}^{\tilde{S}_1}(y+b_1)-\mathbf{\tilde{U}}^*)-(\mathbf{U}^{S_1}(y)-\mathbf{U}_+)]\notag\\
&+[(\mathbf{U}^{\tilde{S}_2}(y+b_2)-\mathbf{U}_+)-(\mathbf{U}^{S_2}(y)-\mathbf{U}^*)], \,\,\,y\geq0.
\end{aligned} \right.
\end{eqnarray}
Obviously, from the Lemma of \ref{key lemma1}, one obtains
\begin{lemma} \label{lemma 3.2}
If the wave strength $\delta$ is suitable small, then for any fixed $b_1,b_2$, it holds that
\begin{equation}\label{lemma 3.2-1}
\begin{split}
&|(\mathbf{U}^{\tilde{S}_1}(y+b_1)-\mathbf{U}_-)-(\mathbf{U}^{S_1}(y)-\mathbf{U}^*)|\leq C_{b_1}\delta^2e^{-\frac12c\delta|y|},\,\,\, \quad  for \,\,y\leq0,\\
&|(\mathbf{U}^{\tilde{S}_2}(y+b_2)-\mathbf{\tilde{U}}^*)-(\mathbf{U}^{S_2}(y)-\mathbf{U}_-)| \leq C_{b_2}\delta^2e^{-\frac12c\delta|y|}, \,\,\,\quad for \,\,y\leq0,\\
&|(\mathbf{U}^{\tilde{S}_1}(y+b_1)-\mathbf{\tilde{U}}^*)-(\mathbf{U}^{S_1}(y)-\mathbf{U}_+)| \leq C_{b_1}\delta^2e^{-\frac12c\delta|y|}, \,\,\,\quad for \,\,y\geq0,\\
&|(\mathbf{U}^{\tilde{S}_2}(y+b_2)-\mathbf{U}_+)-(\mathbf{U}^{S_2}(y)-\mathbf{U}^*)| \leq C_{b_2}\delta^2e^{-\frac12c\delta|y|}, \,\,\,\quad for \,\,y\geq0,
\end{split}
\end{equation}
where the positive constant $C_{b_i}\lesssim e^{|b_i|},i=1,2$.
\end{lemma}
\begin{proof}
Here we only need to give the proof of inequality \eqref{lemma 3.2-1}$_1$,  the proof of the other inequalities are similar. By using (\ref{1-relaxation shock wave}), one has
\begin{eqnarray}
2s_1V^{S_1}_y + (a^2 - s_1^2) V^{S_1}_{yy}=\frac{-s_1^2(V^{S_1}-v^*)-(p(V^{S_1})-p(v^*))}{a^2-s_1^2},
\end{eqnarray}
then, from Lemma \ref{key lemma1}, combining with  \eqref{A change in the intensity of a wave},  one obtains
\begin{eqnarray}\label{v S1 equation}
V^{S_1}_y(y)=\frac{B_1(v^*,v_+,y)}{a^2-s_1^2}(V^{S_1}(y)-v^*),
\end{eqnarray}
with
\begin{eqnarray}
&&B_1(v^*,v_+,y)=\frac1{2|s_1|-O(1)\delta_1}\Big(s_1^2+\frac{p(V^{S_1})-p(v^*)}{V^{S_1}-v^*}\Big)\notag\\
&&=\frac1{2|s_1|-O(1)\delta_1}\Big(\frac{p(V^{S_1})-p(v^*)}{V^{S_1}-v^*}-\frac{p(v_+)-p(v^*)}{v_+-v^*}\Big)\\
&&\approx\frac1{2|s_1|-O(1)\delta_1}\Big(\frac{(\gamma-1)s_1^2}{2v^*}+\frac{\gamma p(v^*)}{{v^*}^2}\Big)[(V^{S_1}-v^*)-(v^+-v^*)]+O(1)\delta^2,\notag
\end{eqnarray}
where  $"\approx"$ means that the term $O(1)\delta B_1$ is omitted, and $c\delta_1\leq B_1\leq C\delta_1$. And then, integrating (\ref{v S1 equation}) with respect to $y$ over $[y,0]$ with $y<0$, one gets
\begin{equation}
v^*-V^{S_1}=(v^*-V^{S_1}(0))e^{-\int_y^0B_1(v^*,v_+,z)dz}=\frac12(v^*-v_+)e^{-\int_y^0B_1(v^*,v_+,z)dz},
\end{equation}
where $V^{S_1}(0)=\frac12(v^*+v_+)$ is used. Through similar derivation,
for $V^{\tilde{S}_1}$, one has
\begin{eqnarray}
v_--V^{\tilde{S}_1}(y+b_1)=\frac12(v^*_--v_-)e^{-\int_y^0\tilde{B}_1(v^*_-,v_-,z)dz},
\end{eqnarray}
with
\begin{eqnarray}
&&\tilde{B}_1(v_-^*,v_-,y)=\frac1{2|\tilde{s}_1|
-O(1)\tilde{\delta}_1}(s_1^2+\frac{p(V^{\tilde{S}_1})-p(v_-)}{V^{\tilde{S}_1}-v_-})\notag\\
&&=\frac1{2|\tilde{s}_1|-O(1)\tilde{\delta}_1}(\frac{p(V^{\tilde{S}_1})-p(v_-)}{V^{\tilde{S}_1}-v_-}-\frac{p(\tilde{v}^*)-p(v_-)}{\tilde{v}^*-v_-})\\
&&\approx\frac1{2|\tilde{s}_1|-O(1)\tilde{\delta}_1}\Big(\frac{(\gamma-1)s_1^2}{2\tilde{v}^*}+\frac{\gamma p(\tilde{v}^*)}{(\tilde{v}^*)^2}\Big)[(V^{\tilde{S}_1}-v_-)-(\tilde{v}^*-v_-)]+O(1)\delta^2,\notag
\end{eqnarray}
and $c\tilde{\delta}_1\leq \tilde{B}_1\leq C\tilde{\delta}_1$. Moreover, one obtains
\begin{eqnarray}\label{prove smallness}
&&|(V^{\tilde{S}_1}(y+b_1)-v_-)-(V^{S_1}(y)-v^*)|\notag\\
&& \leq\frac12\delta_1|e^{-\int_y^0B_1(v^*,v_+,z)dz}-e^{-\int_{y+b_1}^0\tilde{B}_1(v^*_-,v_-,z)dz}|+C\delta^2e^{-c\delta|y|}\notag\\
&& \leq e^{-c\delta|y|}\big(\int_y^0|B_1-\tilde{B}_1|dz+\delta b_1\big)+C\delta^2e^{-c\delta|y|}\notag\\
&& \leq C(1+|b_1|)\Big(\delta^2e^{-c\delta|y|} +\delta e^{-c\delta|y|}\int_y^0|(V^{\tilde{S}_1}-v_-)-(V^{S_1}-v^*)|dz\Big),
\end{eqnarray}
and therefore, applying the Gronwall's inequality to \eqref{prove smallness}, yields that
\begin{eqnarray}\label{inequality v-s}
\int_y^0|(v^{\widetilde{S}_1}-v_-)-(v^{S_1}-v^*)|dz\leq Ce^{|b_1|}\delta,\,\,\,y\leq0.
\end{eqnarray}
Substituting \eqref{inequality v-s} into (\ref{prove smallness}), one has
\begin{eqnarray}
|(V^{\tilde{S}_1}(y+b_1)-v_-)-(V^{S_1}(y)-v^*)|\leq C_{b_1}\delta^2e^{-\frac12c\delta|y|},
\end{eqnarray}
and thus \eqref{lemma 3.2-1} is obtained.  Furthermore, the rest of the inequalities in \eqref{lemma 3.2-1} can be obtained in a similar way,  so far,  the proof of Lemma \ref{lemma 3.2} is completed.
\end{proof}

\

On the basis of  Lemma \ref{lemma 3.2},  one immediately obtains  the following energy perturbation for the shock interaction  $\tilde{\mathbf{U}}^{\tilde{S}_1\tilde{S}_2}-\mathbf{U}^{S_1S_2}$ at the interaction time as follows.
\begin{lemma}\label{pertubation of initail}
It holds that
\begin{eqnarray}\label{key1}
\Big\|\int_{-\infty}^{y}(\tilde{\mathbf{U}}^{\tilde{S}_1\tilde{S}_2}-\mathbf{U}^{S_1S_2})
(z,0)dz\Big\|_{H^2}^2\leq C\delta.
\end{eqnarray}
\end{lemma}

\

On the basis of these preparations, such as  \eqref{shock after t0}, \eqref{initial data after t0} and Lemma 3.3, the problem of (\ref{visous conservation law in section 4}) after the interaction time can be studied in the same way as  section 2, and
the proof is briefly described as follows. Similar as section 2, it is easy to know that, the relaxation shock wave $\tilde{\mathbf{U}}$ after interaction satisfies
\begin{eqnarray}\label{reflected waves}
\tilde{\mathbf{U}}^{\tilde{S}_1\tilde{S}_2}_{\tau}+\mathbf{F}_y(\tilde{\mathbf{U}}^{\tilde{S}_1\tilde{S}_2})=a^2\tilde{\mathbf{U}}^{\tilde{S}_1\tilde{S}_2}_{yy}-\tilde{\mathbf{U}}^{\tilde{S}_1\tilde{S}_2}_{\tau\tau}+\tilde{\mathbf{G}}_y,
\end{eqnarray}
where
\begin{eqnarray}\label{tilde G}
&&\tilde{\mathbf{G}}=\mathbf{F}(\tilde{\mathbf{U}}^{\tilde{S}_1\tilde{S}_2})-\Big(\mathbf{F}(\mathbf{U}^{\tilde{S}_1})+\mathbf{F}(\mathbf{U}^{\tilde{S}_2})-\mathbf{F}(\mathbf{\tilde{U}}^*)\Big),\end{eqnarray}
and
\begin{equation}\label{Inequality of G}
  \tilde{\mathbf{G}}\leq C\delta^3e^{-c\delta|y|-c\delta|\tau|}\overset{\text{def}}{=}O(1)\delta \sigma.
\end{equation}
Setting
\begin{equation}\label{Error in reflected waves}
  \tilde{\mathbf{\Phi}}=\mathbf{U}-\tilde{\mathbf{U}}^{\tilde{S}_1\tilde{S}_2},
\end{equation}
substituting \eqref{Error in reflected waves} into \eqref{visous conservation law in section 4} and (\ref{reflected waves}), one has
\begin{eqnarray}
\mathbf{\tilde\Phi}_{\tau}+\Big(\mathbf{F}(\tilde{\mathbf{U}}^{\tilde{S}_1\tilde{S}_2}+\mathbf{\tilde\Phi})-\mathbf{F}(\tilde{\mathbf{U}}^{\tilde{S}_1\tilde{S}_2})\Big)_y=a^2\mathbf{\tilde\Phi}_{yy}-\tilde{\mathbf{G}}_y+\mathbf{H}_{\tau y}+\mathbf{H}_y-\mathbf{\tilde\Phi}_{\tau\tau}.
\end{eqnarray}
Using the antiderivative method, one defines
\begin{equation}\label{anti-derivative after interfaction}
\tilde{\pmb{\phi}}(y,\tau)=\int_{-\infty}^y\tilde{\mathbf{\Phi}}(z,\tau)dz,
\end{equation}
thus, one obtains
\begin{equation}\label{the system after the wave interaction}
\left\{\begin{split}
&\tilde{\pmb{\phi}}_{\tau}+\mathbf{\tilde{F}}'(\tilde{\mathbf{U}}^{\tilde{S}_1\tilde{S}_2})\tilde{\pmb{\phi}}_y-a^2\tilde{\pmb{\phi}}_{yy}+\tilde{\pmb{\phi}}_{\tau\tau}=
-(\tilde{\mathbf{G}}+\tilde{\mathbf{Q}}),\\
&\chi_\tau=-v\mu,\ \ \ \displaystyle \mu=(\chi^3-\chi)-\big(\frac{\chi_y}{v}\big)_y,
\end{split}\right.\end{equation}
where
\begin{equation}\label{tilde Q}
 \tilde{\mathbf{Q}}=\mathbf{\tilde{F}}(\tilde{\mathbf{U}}^{\tilde{S}_1\tilde{S}_2}+\tilde{\pmb{\phi}}_y)-\mathbf{\tilde{F}}(\tilde{\mathbf{U}}^{\tilde{S}_1\tilde{S}_2})-\mathbf{\tilde{F}}'(\tilde{\mathbf{U}}^{\tilde{S}_1\tilde{S}_2})\tilde{\pmb{\phi}}_y=O(1)|\tilde{\pmb{\phi}}_y|^2.
\end{equation}
Letting
\begin{equation}\label{Diagonalization after the wave interaction}
 \tilde{\pmb{\phi}}={\mathbf{\tilde{R}}}\tilde{\mathbf{w}},
\end{equation}
where
$\tilde{\mathbf{L}}={\mathbf{L}}(\tilde{\mathbf{U}}^{\tilde{S}_1\tilde{S}_2})$, $\tilde{\mathbf{R}}={\mathbf{R}}(\tilde{\mathbf{U}}^{\tilde{S}_1\tilde{S}_2})$ represent the matrix of the left and right eigenvectors respectively, and $\tilde{\mathbf{F}}={\mathbf{F}}(\tilde{\mathbf{U}}^{\tilde{S}_1\tilde{S}_2})$,   diagonalizing the system \eqref{the system after the wave interaction}, one has
\begin{equation}\label{Diagonalized system of equations after interaction}
\left\{
\begin{split}
&\tilde{\mathbf{w}}_{\tau}+\mathbf{\tilde{\Lambda}}\tilde{\mathbf{w}}_y-a^2\tilde{\mathbf{w}}_{yy}+ \tilde{\mathbf{w}}_{\tau\tau}=-\tilde{\mathbf{L}}\big(\tilde{\mathbf{R}}_{\tau}-a^2\tilde{\mathbf{R}}_{yy}+\mathbf{F}'\tilde{\mathbf{R}}_y+\tilde{\mathbf{R}}_{\tau\tau}\big)\tilde{\mathbf{w}}\\
&\qquad\qquad\qquad\qquad\qquad+\tilde{\mathbf{L}}\big(2a^2\tilde{\mathbf{R}}_y\tilde{\mathbf{w}}_y-\tilde{\mathbf{G}}-\tilde{\mathbf{Q}}+\mathbf{H}_{\tau}+\mathbf{H}-2\tilde{\mathbf{R}}_\tau\big)\tilde{\mathbf{w}}_\tau,\\
&\chi_\tau=-v(\chi^3-\chi)+v\big(\frac{\chi_y}{v}\big)_y,\\
&\|\tilde{\mathbf{w}}\|_{H^2}\Big|_{\tau=0}\leq C\delta,\ \ \chi\Big|_{\tau=0}=\chi(y,0),
\end{split}\right.
\end{equation}
where $\tilde{\mathbf{\Lambda}}={\mathbf{\Lambda}}(\tilde{\mathbf{U}}^{\tilde{S}_1\tilde{S}_2})$  is the matrix of eigenvalues for $\tilde{\mathbf{F}}'$.
 Setting the solution space for \eqref{Diagonalized system of equations after interaction} as follows:
\begin{equation}\label{Solution space after interaction}
\begin{split}
& \tilde{X}_{\tilde{M}}(\tilde{I})=\Big\{ (\mathbf{\tilde{w}},\chi)\Big|\mathbf{\tilde{w}}\in C(\tilde{I};H^3), \mathbf{\tilde{w}}_\tau\in C(\tilde{I};H^2), \chi^2-1\in C(\tilde{I};L^2),\chi_y\in C(I;H^3),\\
&\qquad\qquad\qquad\qquad \chi_{y\tau}\in C(\tilde{I};H^2),\mathbf{\tilde{w}}\in L^2(\tilde{I};H^3),\mathbf{\tilde{w}}_\tau\in L^2(\tilde{I};H^2),\\
&\qquad\qquad\qquad\qquad\chi^3-\chi\in L^2(\tilde{I};L^2),\chi_y\in L^2(\tilde{I};H^3),\chi^2-1\in L^2(\tilde{I};L^2),\\
&\qquad\qquad\qquad\qquad\sup_{\tau\in\tilde{I}}\{\|\mathbf{\tilde{w}}\|_3,\|\mathbf{\tilde{w}}_\tau\|_2,\big\|\chi^2-\chi\big\|,\|\chi_y\|_3,\|\chi_{\tau y}\|_1\}\leq M \Big\},
\end{split}\end{equation}
where $\tilde{I}\subseteq [0,+\infty)$ is the any  time interval. Then a priori estimate for the solutions of \eqref{Diagonalized system of equations after interaction} is as following Proposition, the proof is omitted here.
\begin{proposition} \label{a Priori Estimate After the Interaction Time}
There are two constants $\delta_0$ and $M_0$ such that, if $\delta\leq\delta_0$, $M<M_0$,  there exists
a unique strong solution $\mathbf{w}$ to (\ref{Diagonalized system of equations after interaction}) in $\tilde{X}_{\tilde{M}}(0,\infty)$. Furthermore, it holds that
\begin{equation}\label{Energy estimates after the collision}
\begin{split}
&\|\mathbf{\tilde{w}}(\tau)\|_{H^3}^2+\|\mathbf{\tilde{w}}_\tau(\tau)\|_{H^2}^2+\|\chi^2(\tau)-1\|^2+\|\chi_{y}(\tau)\|_3^2+\|\chi_{y\tau}(\tau)\|_2^2\\
&+\int_0^{\infty}\big(\|\mathbf{w}_y(\tau)\|_{H^2}^2+\|\mathbf{w}_\tau(\tau)\|_{H^2}^2+\|(\chi^3-\chi)(\tau)\|^2+\|\chi_y(\tau)\|_3^2\big)d\tau\\
&\leq CM_0(\delta+1),
\end{split}
\end{equation}where $C$ is the positive constants which may depend on $(v_-,u_-)$.
\end{proposition}

\

By using $\pmb{\phi}=\mathbf{R}(\mathbf{U}^{S_1S_2})\mathbf{w}$ before the interaction time $t_0$,  and $\tilde{\pmb{\phi}}=\mathbf{R}(\tilde{\mathbf{U}})\tilde{\mathbf{w}}$ after the interaction time $t_0$, combining with Propositions 2.1 and 3.1,  Theorem 1.1. is immediately obtained, the details are omitted.

\section{Numerical Verification}\label{sec3}
\begin{figure}[htbp]
		\centering
		\includegraphics[width=2.8in]{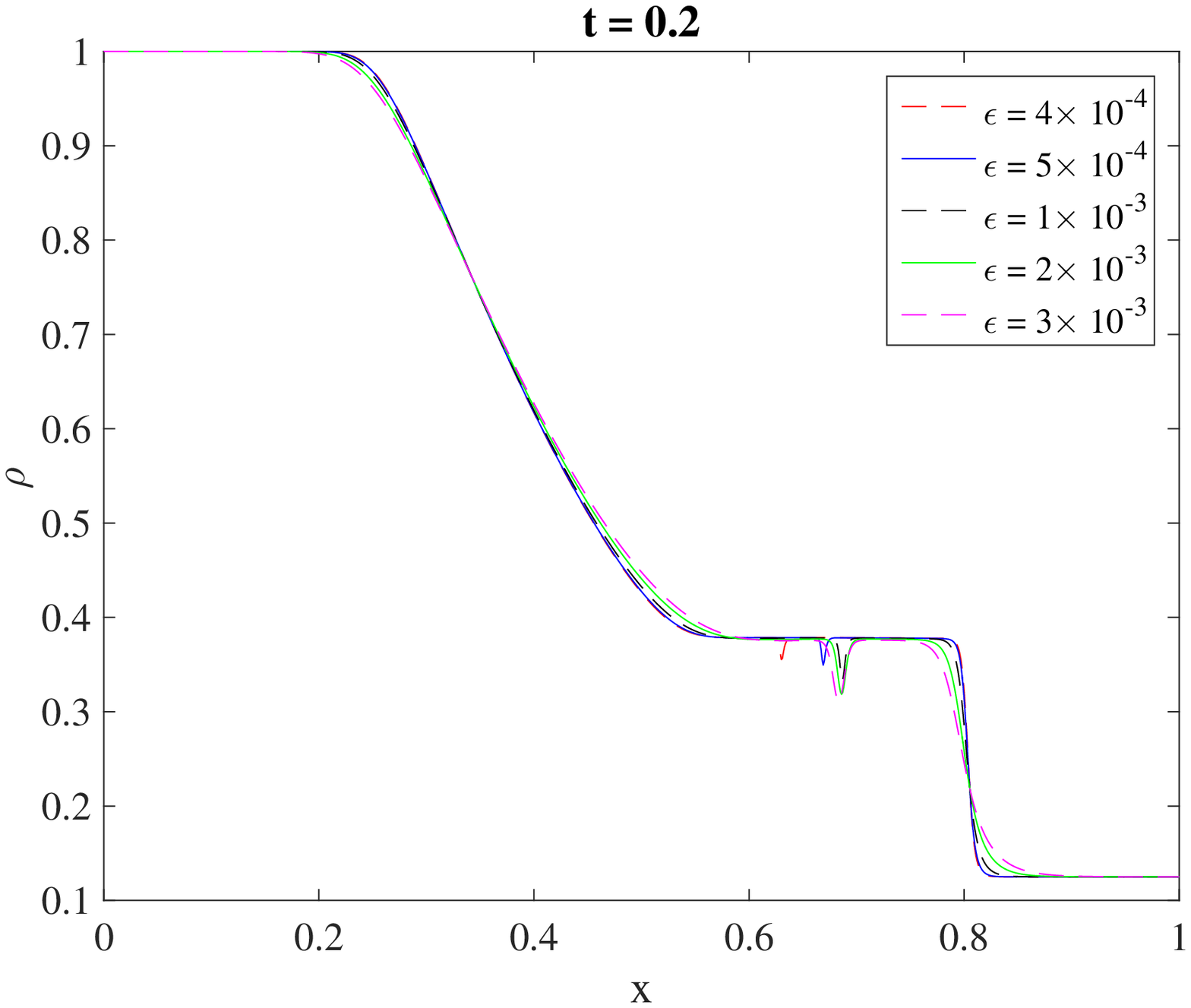}\quad \quad
	\includegraphics[width=2.8in]{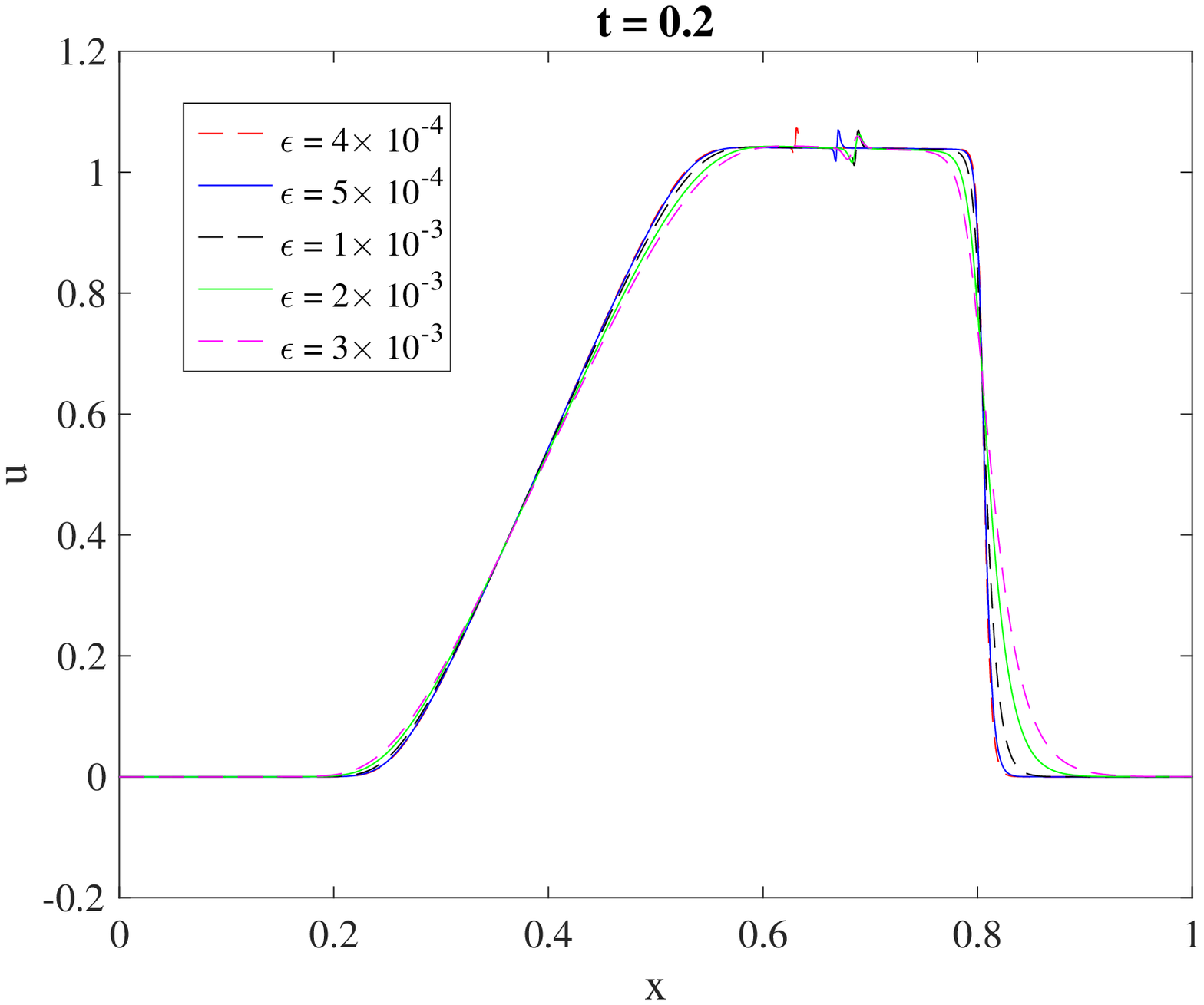}\\
		\includegraphics[width=2.8in]{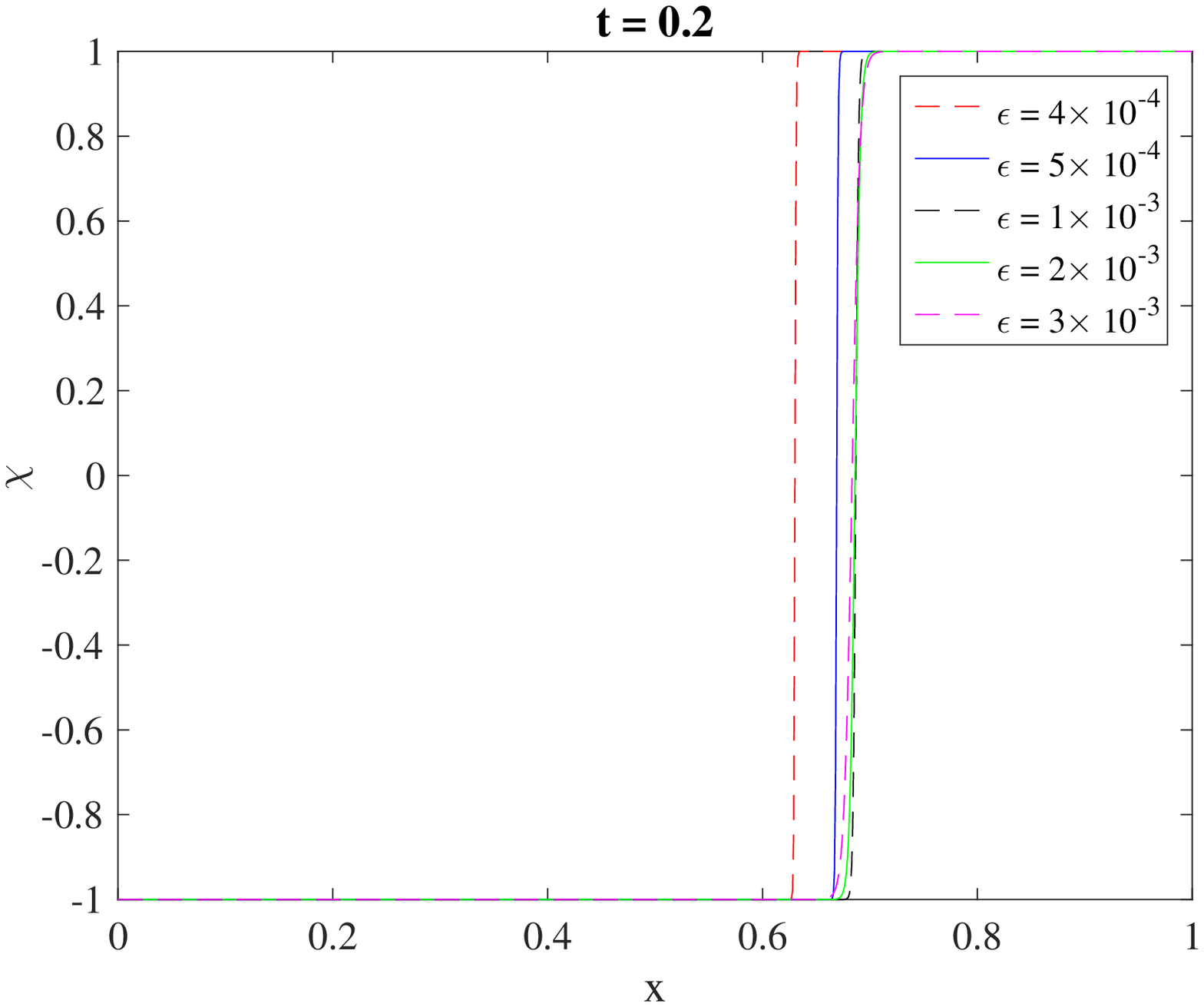}
		\caption{Plot of  $\rho$, $u$ and $\chi$ at time $t = 0.2$. }\label{fig1}
\end{figure}

\setcounter{equation}{0}

\ \ \ \ In this section, we present a numerical example to study  the  relation between shock wave solution and the sharp interface limit for the Navier-Stokes/Allen-Cahn system \eqref{original NSAC} in 1-D.  We use the Jin-Xin relaxation scheme (see Jin-Xin \cite{JinXin1995}) for the first two equations in \eqref{original NSAC} and stabilized method for the third equation in \eqref{original NSAC}. This composite numerical schemes is first proposed and analyzed by He-Shi \cite{NSCHnume}, also see \cite{H-L-S-2018}, \cite{NSACnume}. 
The following example is given by using this composite numerical method.  The computational domain is $[0,1]$ and the initial data is
\begin{align*}
&  (\rho_{0}, u_{0}, \chi_{0})^{T} = \left(1.0, 0.0, \tanh\left(\frac{x-0.5}{0.1\sqrt{5.0}}\right)\right)^{T},   \  0 \leq x < 0.5,\\
	&  (\rho_{0}, u_{0}, \chi_{0})^{T} = \left(0.125, 0.0, \tanh\left(\frac{x-0.5}{0.1\sqrt{5.0}}\right)\right)^{T},   \  0.5 \leq x \leq 1.0.
\end{align*}
We take spatial interval $\Delta x = 1/1000$ and the  parameter $\epsilon $ varies from $4\times 10^{-4}$ to $3 \times 10^{-3}$. Figure \ref{fig1} show the numerical solutions of $\rho$, $u$ and $\chi$  when $t = 0.2$. It is obvious to observe   that the interface  does not affect the shock wave solution.  When $\epsilon$ becomes smaller,  the shock wave solution becomes more sharp  due to the kinematic viscosity coefficient $\nu(\epsilon)$.  When $\epsilon$ changes,  the shock wave solution is developed  meanwhile the sharp interface limit  is approached.

\begin{figure}[htbp]
	\centering
	\includegraphics[width=2.8in]{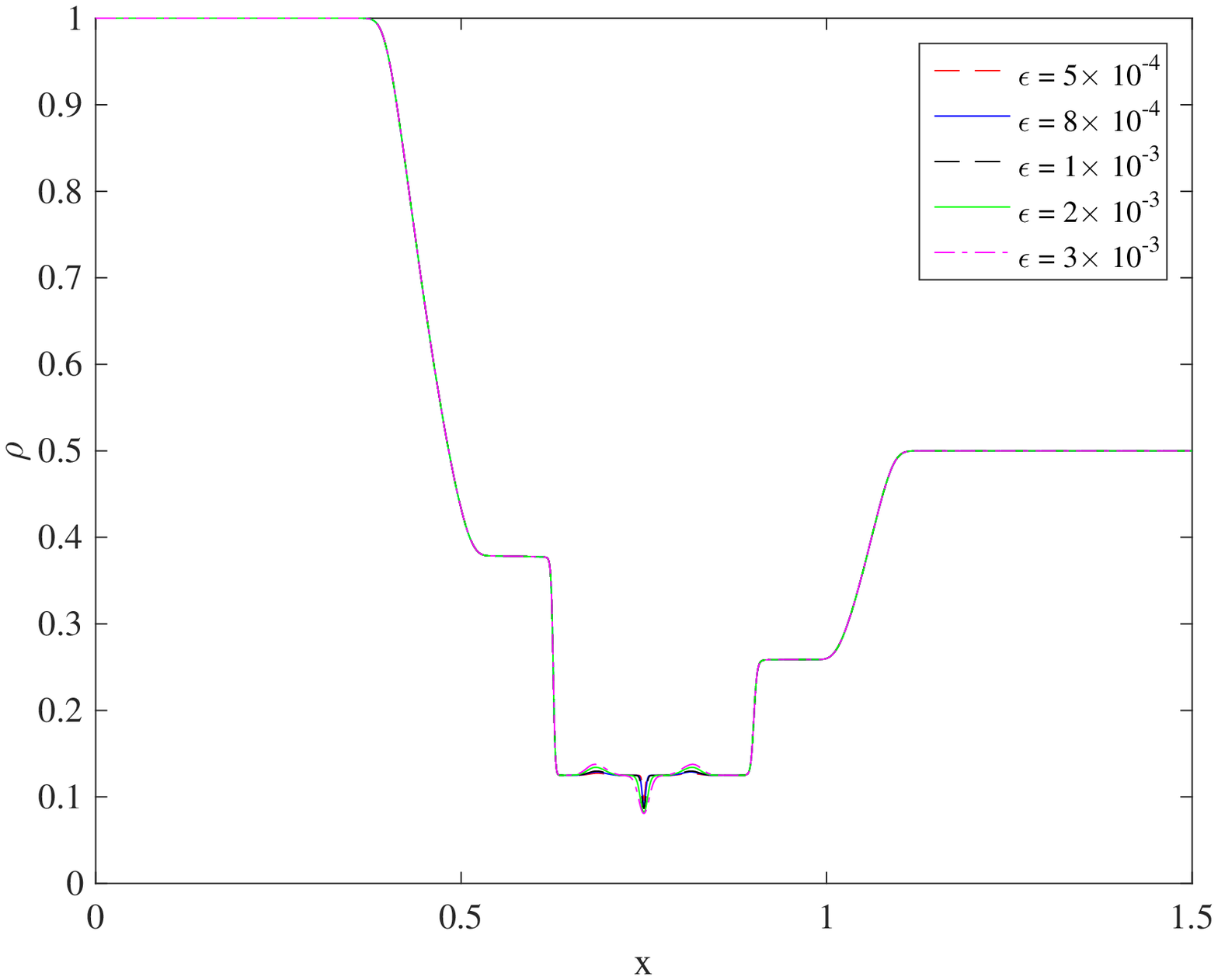}\quad \quad
	\includegraphics[width=2.8in]{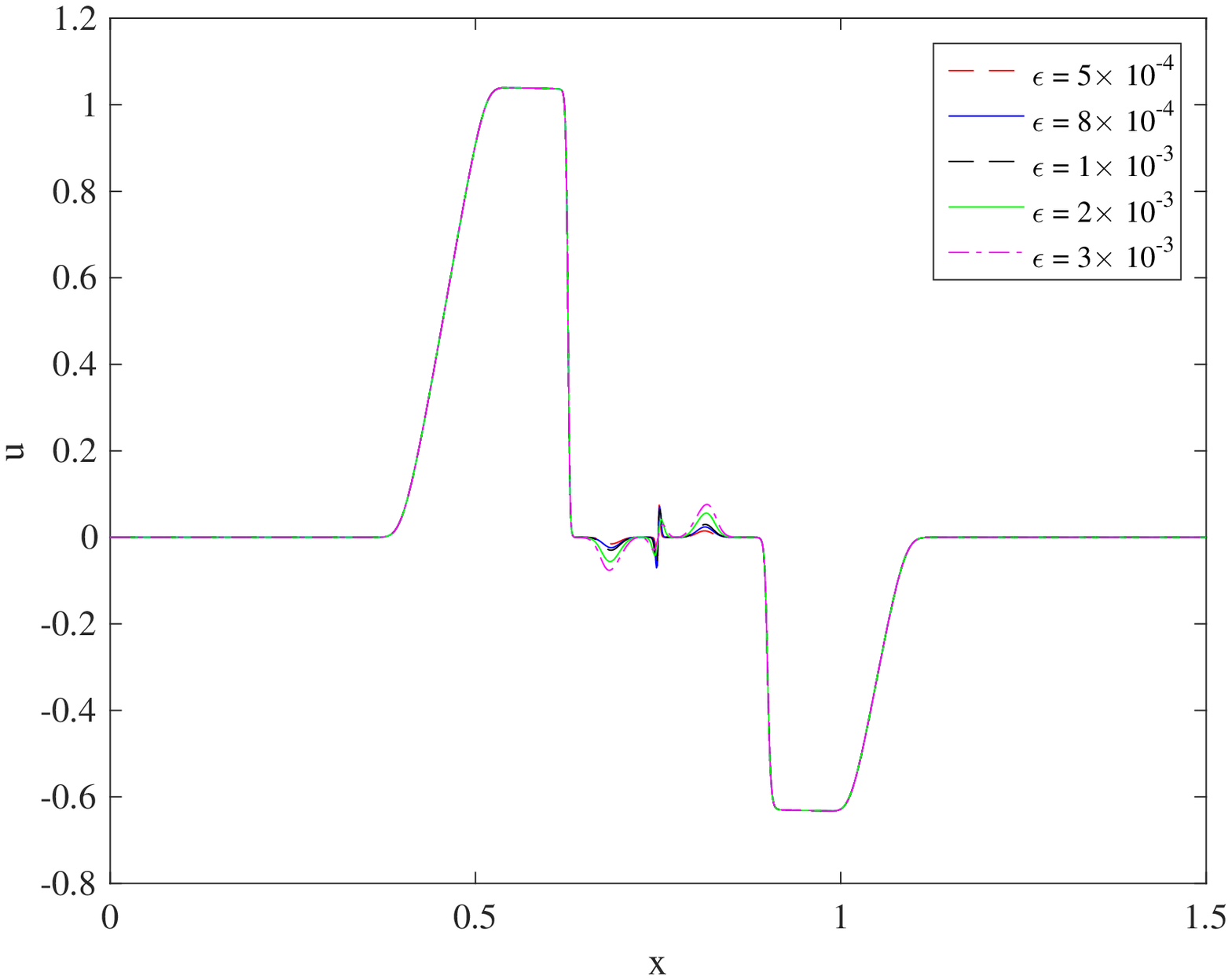}\\
	\includegraphics[width=2.8in]{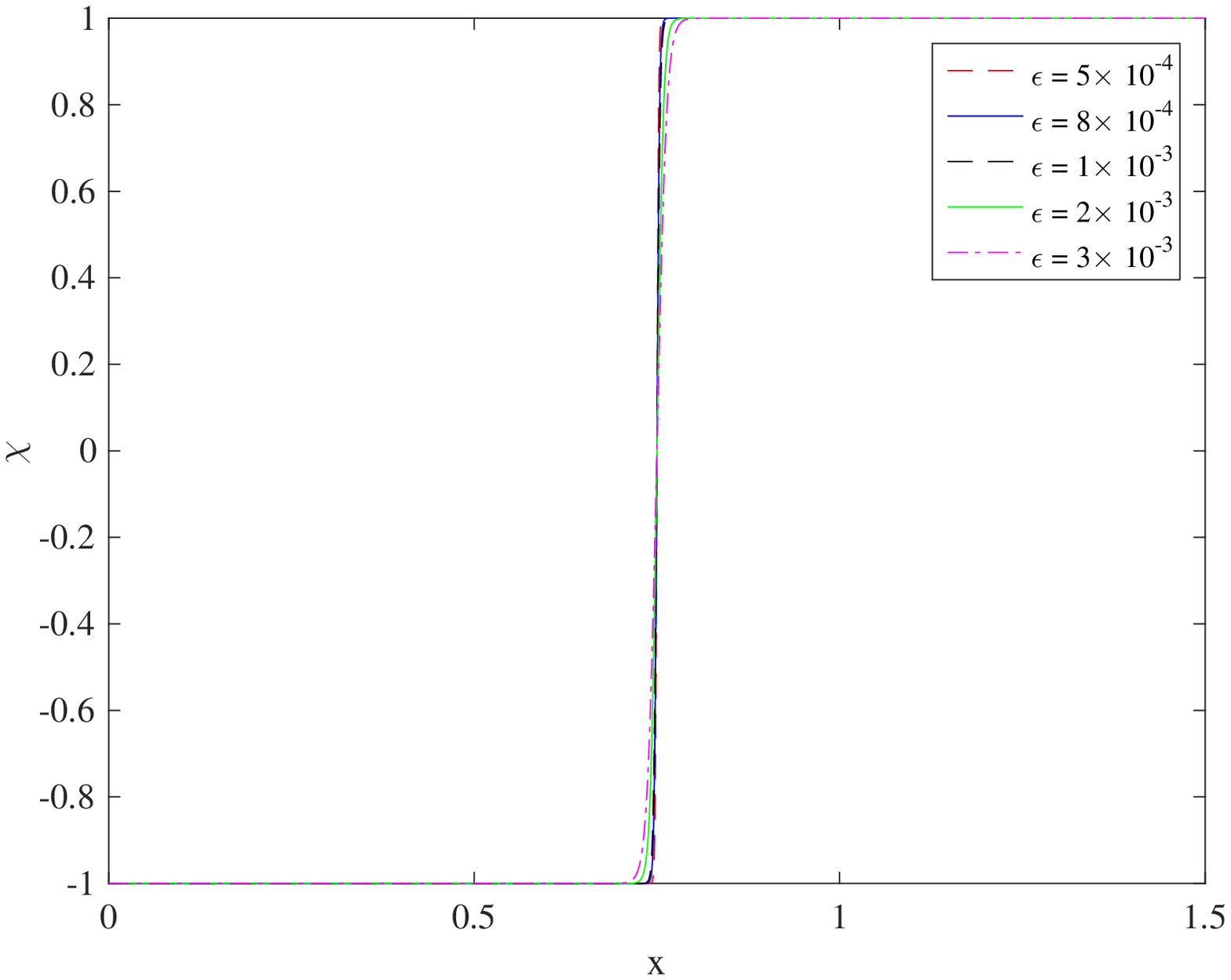}
	\caption{Plot of  $\rho$, $u$ and $\chi$ at time $t = 0.08$ (before interaction). }\label{fig2}
\end{figure}

\begin{figure}[htbp]
	\centering
	\includegraphics[width=2.8in]{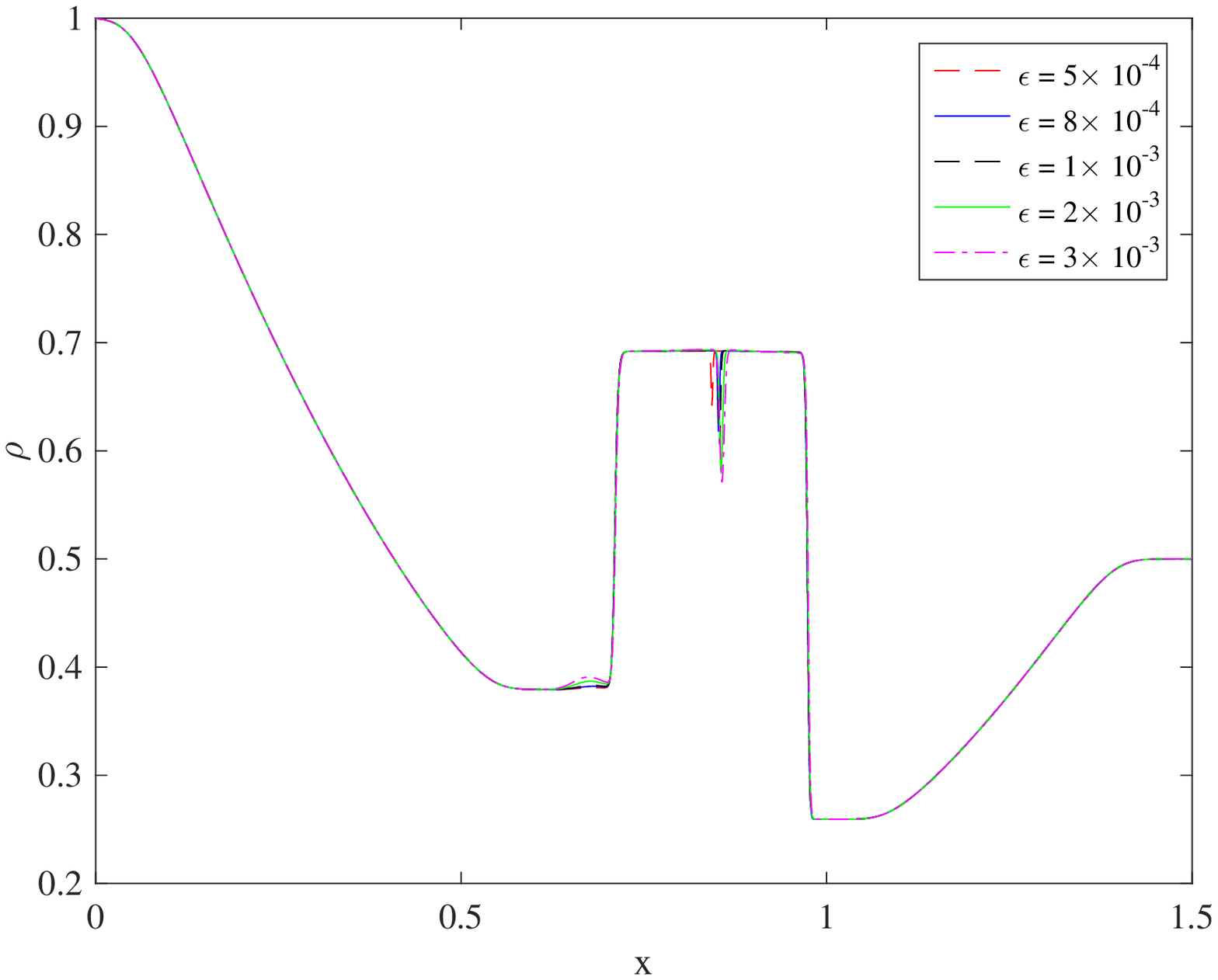}\quad \quad
	\includegraphics[width=2.8in]{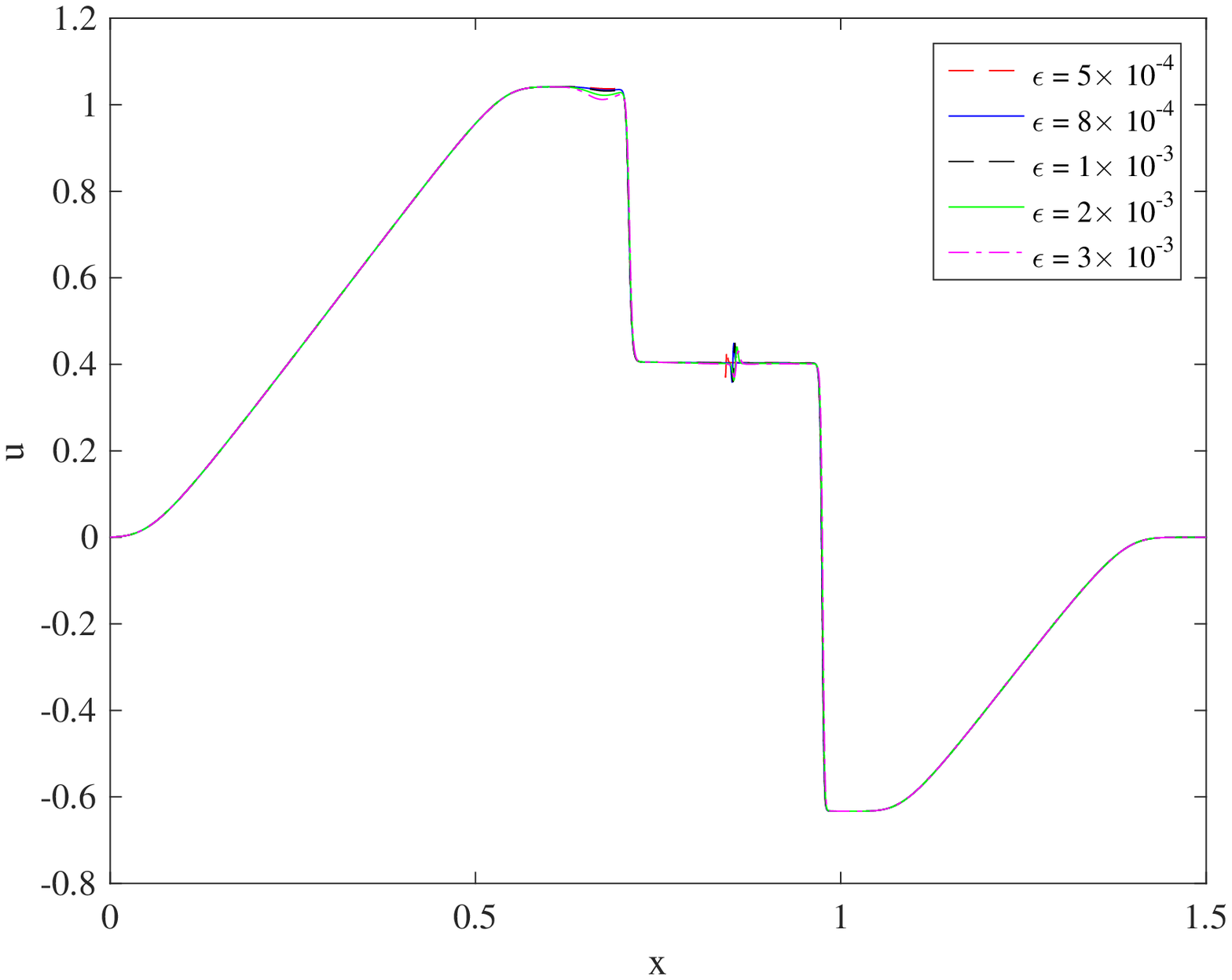}\\
	\includegraphics[width=2.8in]{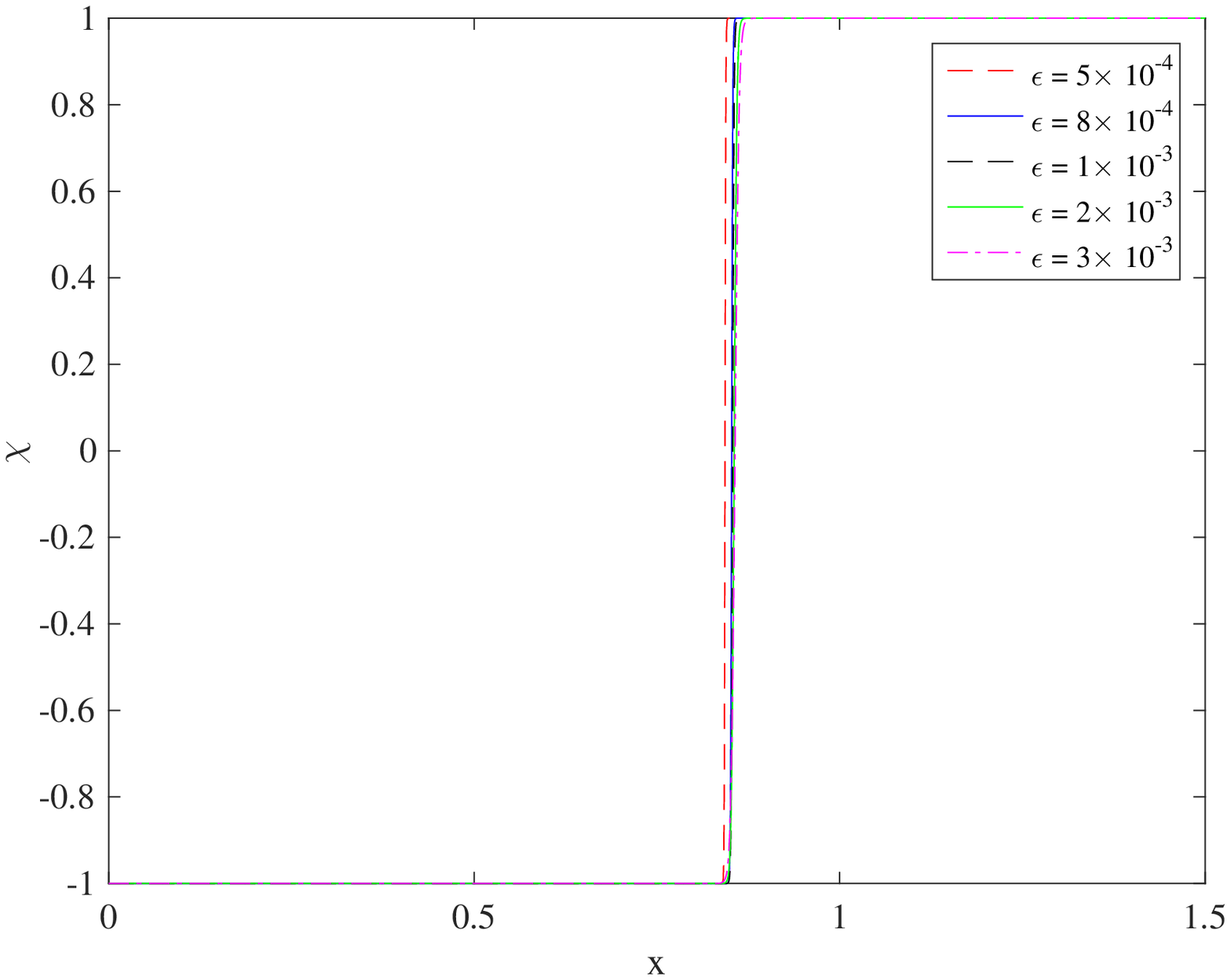}
	\caption{Plot of  $\rho$, $u$ and $\chi$ at time $t = 0.4$ (after interaction). }\label{fig3}
\end{figure}

We also consider two waves interaction case. The computational domain is $[0,1.5]$ and the initial data is
\begin{align*}
&  (\rho_{0}, u_{0}, \chi_{0})^{T} = \left(1.0, 0.0, \tanh\left(\frac{x-0.75}{0.1\sqrt{5.0}}\right)\right)^{T},   \  0 \leq x < 0.5,\\
&  (\rho_{0}, u_{0}, \chi_{0})^{T} = \left(0.125, 0.0, \tanh\left(\frac{x-0.75}{0.1\sqrt{5.0}}\right)\right)^{T},   \  0.5 \leq x < 1.0,\\
& (\rho_{0}, u_{0}, \chi_{0})^{T} = \left(0.5, 0.0, \tanh\left(\frac{x-0.75}{0.1\sqrt{5.0}}\right)\right)^{T},   \  1.0 \leq x \leq 1.5.
\end{align*}
We take spatial interval $\Delta x = 1.5/2000$ and the  parameter $\epsilon $ varies from $5\times 10^{-4}$ to $3 \times 10^{-3}$.  Figure \ref{fig2} show the numerical results of $\rho$, $u$ and $\chi$  when $t = 0.08$,  which is the case before interaction. Figure \ref{fig3} show the numerical results of $\rho$, $u$ and $\chi$  when $t = 0.4$, which is the case after interaction. Therefore, the interaction of two shock does not affect the interface.


\end{document}